# EFFICIENT LINEAR, STABILIZED, SECOND-ORDER TIME MARCHING SCHEMES FOR AN ANISOTROPIC PHASE FIELD DENDRITIC CRYSTAL GROWTH MODEL

XIAOFENG YANG

ABSTRACT. We consider numerical approximations for a phase field dendritic crystal growth model, which is a highly nonlinear system that couples the anisotropic Allen-Cahn type equation and the heat equation together. We propose two efficient, linear, second-order time marching schemes. The first one is based on the linear stabilization approach where all nonlinear terms are treated explicitly and one only needs to solve two linear and decoupled second-order equations. Two linear stabilizers are added to enhance the energy stability, therefore the scheme is quite efficient and stable that allows for large time steps in computations. The second one combines the recently developed Invariant Energy Quadratization approach with the linear stabilization approach. Two linear stabilization terms, which are shown to be crucial to remove the oscillations caused by the anisotropic coefficients numerically, are added as well. We further show the obtained linear system is well-posed and prove its unconditional energy stability rigorously. For both schemes, various 2D and 3D numerical simulations are implemented to demonstrate the stability and accuracy.

## 1. INTRODUCTION

A dendritic crystal that develops with a typical multi-branching tree-like form, is very common in daily life, for instances, snowflake formation, frost patterns on a substrate, solidification of metals, or crystallization in supersaturated solutions, etc., see [5, 8, 10]. In the modeling and numerical simulations for investigating the process of dendritic crystal growth, the phase field method, as a powerful tool for simulating free interfacial motions, had been widely used, see the pioneering modeling work of Halperin, Kobayashi, and Collins et. al. in [3, 6, 15], and a number of subsequent modeling/simulation advances in [1, 12–14, 18, 20, 23, 29, 40].

In this paper, we consider numerical approximations for a phase field dendritic crystal growth model which was proposed in [14]. The system consists of two coupled, second-order equations: the Allen-Cahn type equation with a gradient-dependent anisotropic coefficient, and the heat transfer equation. More precisely, an order parameter (phase field variable) is introduced to define the physical state (liquid or solid) of the system at each point, and a free-energy functional is devised by incorporating a specific form of the conformational entropy with anisotropic spatial gradients that is introduced to study dynamics of atomic-scale dendritic crystal growth. The evolution equation for the phase field variable follows the $L^2$-gradient flow of the total free energy, and the coupling between the Allen-Cahn equation and heat equation is realized by a nonlinear heat-transfer term along the interface. Comparing to the isotropic phase field model, in addition to the stiffness issue induced by the interfacial width, there are a couple more challenges for solving the anisotropic dendritic crystal system, including (i) how to discretize the anisotropic coefficient which can induce large oscillations (cf. Fig. 4.2); and (ii) how to discretize the nonlinear coupling heat transfer term. These difficulties increase the complexity of algorithm developments to a large extent.

Therefore, there are very few successful attempts in designing efficient and energy stable schemes because of the complexity of the nonlinear terms in the anisotropic dendritic model. Most of the available schemes are either nonlinear which need some efficient iterative solvers, and/or do not preserve energy stability at all (cf. [1, 11, 16, 17, 19, 21, 22] and the references therein). In [40], the authors applied the Invariant Energy Quadratization (IEQ) approach to develop a second-order, energy stable scheme by enforcing the free-energy density as an invariant, quadratic functional in terms of new, auxiliary variables and treating the nonlinear terms semi-explicitly. However, the energy stability preserves only when the temperature is a constant. Meanwhile, even though the discrete energy law holds formally, it is hard to show the discrete energy obtained in [40] is bounded from below.

In this paper, we aim to develop some more effective and efficient numerical schemes to solve the dendritic phase field model [14]. We expect that our schemes can combine the following three desired properties, i.e., (i) accurate (second order in time); (ii) energy stable (the discrete energy must be bounded from below while it follows the unconditional energy law); and (iii) easy to implement (linear). To this end, we develop two linear, second-order time marching schemes. The first one is based on the stabilized-explicit (or called linear stabilization) approach (cf. [9, 25–27, 30, 31]), where all nonlinear terms are treated explicitly and one only needs to solve two decoupled, second-order equations. To enhance the energy stability, two linear stabilizers are added, therefore the scheme is quite efficient and stable that allows for the large time step in computations. However, its energy stability is not provable yet. The second scheme is based on the stabilized Invariant Energy Quadratization (SIEQ) approach which combines the recently developed IEQ approach (cf. [2, 7, 32–36, 38–40]) with the linear stabilization approach. That is to say, besides the quadratization technique of using the auxiliary variable, the two identical linear stabilizers are added as well. Numerical examples show that these two stabilizing terms can efficiently eliminate oscillations and thus allow large time steps. We rigorously prove that the SIEQ scheme is unconditionally energy stable and the obtained linear system is well-posed. To the best of the author's knowledge, the developed SIEQ scheme is the first scheme for the anisotropic dendritic system that can be theoretically proved to be energy stable. Through the comparisons with two other prevalent numerical schemes such as the fully-explicit and the IEQ schemes (without stabilizers) for a number of classical benchmark numerical examples, we demonstrate the stability and the accuracy of the proposed schemes as well.

The rest of the paper is organized as follows. In Section 2, we give a brief introduction of the governing PDE system for the anisotropic dendritic model. In Section 3, we develop two numerical schemes with the second-order accuracy for simulating the model, and rigorously prove the well-posedness and the unconditional energy stabilities for the second one. Various 2D and 3D numerical experiments are given in Section 4 to validate the accuracy and efficiency of the proposed numerical schemes. Finally, some concluding remarks are given in Section 5.

## 2. Model equations and its energy law

Now we give a brief description for the anisotropic phase field dendritic crystal growth model proposed in [14]. Let $\Omega$ be a smooth, open, bounded, connected domain in $\mathbb{R}^d$ with $d = 2, 3$. We introduce a scalar function $\phi(\boldsymbol{x}, t)$ which is an order parameter to label the liquid and solid phase, where $\phi = 1$ for the solid and $\phi = -1$ for the fluid. These two values are connected by a smooth transitional layer with the thickness $\epsilon$.

We consider the total free energy as follows,

$$(2.1) \qquad E(\phi, u) = \int_\Omega \Big(\frac{1}{2}|\kappa(\nabla\phi)\nabla\phi|^2 + \frac{F(\phi)}{4\epsilon^2} + \frac{\lambda}{2\epsilon K}u^2\Big)d\boldsymbol{x},$$

in which, $\epsilon$, $\lambda$ and $K$ are all positive parameters, $u(\boldsymbol{x}, t)$ is the temperature, $F(\phi) = (\phi^2 - 1)^2$ is the double well potential, $\kappa(\nabla\phi)$ is a function describing the anisotropic property that depends on on the direction of the outer normal vector $\boldsymbol{n}$ which is the interface normal defined as $\boldsymbol{n} = -\frac{\nabla\phi}{|\nabla\phi|}$. For the 2D system, the anisotropy coefficient $\kappa(\nabla\phi)$ is usually given by

$$(2.2) \qquad \kappa(\nabla\phi) = 1 + \epsilon_4 \cos(m\Theta),$$

where $m$ is the number of folds of anisotropy, $\epsilon_4$ is the parameter for the anisotropy strength, and $\Theta = \arctan(\frac{\phi_y}{\phi_x})$. when $m = 4$ (i.e., fourfold anisotropy), for instance, $\kappa(\nabla\phi)$ can be easily reformulated in terms of the phase field variable $\phi$, namely, for 2D,

$$(2.3) \qquad \kappa(\nabla\phi) = (1 - 3\epsilon_4)\Big(1 + \frac{4\epsilon_4}{1 - 3\epsilon_4}\frac{\phi_x^4 + \phi_y^4}{|\nabla\phi|^4}\Big);$$

and for 3D,

$$(2.4) \qquad \kappa(\nabla\phi) = (1 - 3\epsilon_4)\Big(1 + \frac{4\epsilon_4}{1 - 3\epsilon_4}\frac{\phi_x^4 + \phi_y^4 + \phi_z^4}{|\nabla\phi|^4}\Big).$$

By adopting the Allen-Cahn type ($L^2-$gradient flow) relaxation dynamics for the dendritic crystal growth, one obtains the governing dynamical equations via the variational approach, which reads as follows:

$$(2.5) \qquad \tau(\phi)\phi_t = -\frac{\delta E}{\delta \phi} - \frac{\lambda}{\epsilon}p'(\phi)u,$$

$$= \nabla \cdot (\kappa^2(\nabla\phi)\nabla\phi + \kappa(\nabla\phi)|\nabla\phi|^2 \boldsymbol{H}(\phi)) - \frac{f(\phi)}{\epsilon^2} - \frac{\lambda}{\epsilon}p'(\phi)u,$$

$$(2.6) \qquad u_t = D\Delta u + Kp'(\phi)\phi_t,$$

where $\boldsymbol{H}(\phi)$ is the variational derivative of $\kappa(\nabla\phi)$. In 2D, it reads as

$$(2.7) \qquad \boldsymbol{H}(\phi) = \frac{\delta\kappa(\nabla\phi)}{\delta\phi} = 4\varepsilon_4\frac{4}{|\nabla\phi|^6}\Big(\phi_x(\phi_x^2\phi_y^2 - \phi_y^4), \phi_y(\phi_x^2\phi_y^2 - \phi_x^4)\Big),$$

whereas in 3D,

$$\boldsymbol{H}(\phi) = \frac{\delta \kappa(\nabla \phi)}{\delta \phi} = 4\varepsilon_4 \frac{4}{|\nabla \phi|^6} \Big(\phi_x(\phi_x^2 \phi_y^2 + \phi_x^2 \phi_z^2 - \phi_y^4 - \phi_z^4),$$
$$\phi_y(\phi_y^2 \phi_z^2 + \phi_x^2 \phi_y^2 - \phi_x^4 - \phi_z^4), \ \phi_z(\phi_x^2 \phi_z^2 + \phi_y^2 \phi_z^2 - \phi_x^4 - \phi_y^4)\Big). \quad (2.8)$$

In the system (2.5)-(2.6), $\tau(\phi) > 0$ is the mobility constant that is chosen either as a constant [14,29], or as a function of $\phi$ [14], $D$ is the diffusion rate of the temperature, $\frac{\delta E}{\delta \phi}$ is the variational derivative of the total energy with respect to $\phi$, the function $p(\phi)$ that accounts for the generation of latent heat that is a phenomenological functional taking the form that it preserves the minima of $\phi$ at $\pm 1$ independently of the local value of $u$. For $p(\phi)$, there are two common choices: $p(\phi) = \frac{1}{5}\phi^5 - \frac{2}{3}\phi^3 + \phi$ and $p'(\phi) = (1 - \phi^2)^2$ (cf. [14,15]); or $p(\phi) = \phi - \frac{1}{3}\phi^3$ and $p'(\phi) = 1 - \phi^2$ (cf. [28]).

Without the loss of generality, we adopt the periodic boundary condition or the no-flux homogenous Neumann boundary conditions in order to remove all complexities associated with the boundary integrals in this study, i.e.,

$$(2.9) \qquad \text{(i) all variables are periodic; or (ii)} \ \frac{\partial \phi}{\partial \mathbf{n}}\Big|_{\partial \Omega} = \frac{\partial u}{\partial \mathbf{n}}\Big|_{\partial \Omega} = 0,$$

where $\mathbf{n}$ is the outward normal of the computational domain $\Omega$.

The model equations (2.5)-(2.6) follows the dissipative energy law. By taking the $L^2$ inner product of (2.5) with $\phi_t$, and of (2.6) with $\frac{\lambda}{\epsilon K}T$, using the integration by parts and combining the obtained two equalities, we obtain

$$(2.10) \qquad \frac{d}{dt}E(\phi, u) = -\|\sqrt{\tau(\phi)}\phi_t\|^2 - \frac{\lambda D}{\epsilon K}\|\nabla u\|^2 \leq 0.$$

**Remark 2.1.** *Note in many literatures, $p'(\phi)$ in (2.5) is chosen to be $(1 - \phi^2)^2$ or $1 - \phi^2$, but $p'(\phi)$ in (2.6) is set to be 1, instead. At this time, the whole system does not preserve the energy dissipation law any more, see the detailed discussions in [14].*

3. NUMERICAL SCHEMES

Now we aim to develop efficient schemes for solving the dendritic model (2.5)-(2.6) where the main challenging issues are to find proper approaches to discretize nonlinear terms, particularly, the terms associated with the anisotropic gradient-dependent coefficient $\kappa(\nabla \phi)$, as well as the nonlinear coupling term between the phase field variable and temperature ($Tp'(\phi)$ in (2.5) and $p'(\phi)\phi_t$ in (2.6)).

We recall that there exist quite a few numerical techniques that can preserve the unconditional energy stability, for the isotropic models, e.g., the isotropic Allen-Cahn (or Cahn-Hilliard) equation where the only numerical challenge is to discretize the cubic term from the double-well potential, for instances, the convex splitting approach [4, 24], the linear stabilized-explicit approach [9, 25–27, 30, 31], and the recently developed IEQ approach [2, 7, 32–36, 38–40], etc. In the convex splitting approach, the convex part of the energy potential is treated implicitly and the concave part is treated explicitly. While it is provably energy stable, it usually produces a nonlinear scheme in most cases, thus the implementation is complicated and the computational cost is high. Meanwhile, for the specific dendritic model, it is questionable whether the gradient potential multiplied with the anisotropic coefficient could be split into the combinations of the convex and concave parts.

In this paper, we focus on developing *linear* schemes that possess advantages of easy implementations and lower computational cost. The first proposed scheme is based on the fully-explicit approach, where all nonlinear terms are treated explicitly. However, it is well-known that fully-explicit type schemes can lead to considerable time step constraints, see the numerical examples in Section 4.1 and 4.2, as well as more theoretical/numerical evidence given in [26] for isotropic phase field models. For the anisotropic dendritic model, things are about to get worse since the anisotropic coefficient can still cause oscillations even when using very small time steps (shown in Fig.4.4(b)). Nonetheless, we must not lose sight of the fact that the fully-explicit schemes had been widely used for solving the dendritic phase field model, see [1, 12–14, 18, 20, 23, 29]. Therefore, to fix the inherent deficiency of the fully-explicit schemes, we propose the stabilized-explicit scheme where two linear stabilizers are added into the explicit scheme and present a number of numerical simulations to see how the stability gets improved. While numerical examples show that this proposed stabilized-explicit scheme is quite stable for all trial time steps, it is still an open question on how to prove the energy stability theoretically when the anisotropic coefficient is involved.

The second proposed scheme is based on the recently developed IEQ approach. Similar to the explicit type scheme, even though the IEQ scheme is formally unconditionally energy stable for the anisotropic case, it still blows up when large time steps are used due to the large spacial oscillations caused by anisotropy (cf. Fig. 4.1 and Fig. 4.5). Therefore, we propose the stabilized-IEQ approach where the IEQ approach is combined with the linear stabilization technique. By adding two linear stabilizing terms and treating all involved nonlinear terms in the *semi-explicit* way, a well-posed linear system is obtained and it can be proved that the system is unconditionally energy stable. These two stabilizing terms are not only crucial to enhance energy stability in computations, but also are the keys to prove the well-posedness of the linear system.

3.1. **Stabilized-Explicit scheme.** We first present a second-order, stabilized-explicit time discretization scheme, as follows.

**Scheme 1.** *Assuming $\phi^n$ and $u^n$ are known, we update $\phi^{n+1}$ and $u^{n+1}$ from the following two steps:*

**step 1:**

$$
\begin{aligned}
\tau(\phi^{\star,n+1})\frac{3\phi^{n+1} - 4\phi^n + \phi^{n-1}}{2\delta t} &+ \frac{S_1}{\epsilon^2}(\phi^{n+1} - 2\phi^n + \phi^{n-1}) - S_2\Delta(\phi^{n+1} - 2\phi^n + \phi^{n-1}) \\
= \nabla \cdot &\left( \kappa^2(\nabla\phi^{\star,n+1})\nabla\phi^{\star,n+1} + \kappa(\nabla\phi^{\star,n+1})|\nabla\phi^{\star,n+1}|^2 \boldsymbol{H}^{\star,n+1} \right) \\
&- \frac{f(\phi^{\star,n+1})}{\epsilon^2} - \frac{\lambda}{\epsilon}p'(\phi^{\star,n+1})u^{\star,n+1};
\end{aligned}
$$

(3.1)

**step 2:**

$$
\frac{3u^{n+1} - 4u^n + u^{n-1}}{2\delta t} - D\Delta u^{n+1} = Kp'(\phi^{n+1})\frac{3\phi^{n+1} - 4\phi^n + \phi^{n-1}}{2\delta t}.
$$

(3.2)

Here, $S_i$, $i = 1, 2$ are two positive stabilizing parameters, and

$$(3.3) \qquad \phi^{\star,n+1} = 2\phi^n - \phi^{n-1}, \; u^{\star,n+1} = 2u^n - u^{n-1}, \; \boldsymbol{H}^{\star,n+1} = \boldsymbol{H}(\phi^{\star,n+1}).$$

In the above scheme, all nonlinear terms are treated explicitly. To enhance the stability, we add two second-order linear stabilizers (associated with $S_1$ and $S_2$) in (3.1). When both of them vanish, the scheme becomes fully-explicit type scheme that had been widely used in literatures, see [1, 12–14, 18, 20, 23, 29].

The first stabilizing term, $\frac{S_1}{\epsilon^2}(\phi^{n+1} - 2\phi^n + \phi^{n-1})$, is actually a well-known trick used in the linear stabilization approach to balance the explicit treatment of the $\frac{1}{\epsilon^2} f(\phi)$ for solving the isotropic phase field model (cf. [26]). Similarly, the second stabilizer term, $-S_2 \Delta(\phi^{n+1} - 2\phi^n + \phi^{n-1})$, is used to balance the explicit treatment of the gradient term with the anisotropic coefficient. The errors that these two terms introduced are of order $\frac{S_1}{\epsilon^2} \delta t^2 \phi_{tt}(\cdot)$ and $S_2 \delta t^2 \Delta \phi_{tt}(\cdot)$, respectively, which are of the same order as the error introduced by the second order extrapolation for the nonlinear term $f(\phi)$ and the anisotropic gradient term. Numerical examples show that these two stabilizers are crucial to removing all oscillations induced by the anisotropic coefficient $\gamma(\boldsymbol{n})$, since the term $\boldsymbol{n} = -\frac{\nabla \phi}{|\nabla \phi|}$ changes its sign frequently in the bulk part where $|\nabla \phi|$ is close to zero (cf. Fig. 4.2).

As mentioned above, the stabilized explicit scheme is quite efficient and one only needs to solve two decoupled second-order type equations, however, we are not aware how to prove its energy stability theoretically even though its computational performance is almost as good as the stabilized-IEQ scheme with the provable energy stability (cf. numerical examples 4.1 and 4.2).

3.2. **Stabilized-IEQ scheme.** The key procedure of the IEQ approach is to qudratize the free energy potential via one or more auxiliary variables. Thus we define an auxiliary variable as follows:

$$(3.4) \qquad U = \sqrt{\frac{1}{2}|\kappa(\nabla \phi)\nabla \phi|^2 + \frac{1}{4\epsilon^2} F(\phi) + B},$$

where $B$ is a constant that can ensure the radicand positive (in all numerical examples, we let $B = 5 \times 10^4$ which is the same order of $\frac{1}{\epsilon^2}$). Thus the total free energy (2.1) can be rewritten as

$$(3.5) \qquad E(\phi, U, u) = \int_\Omega (U^2 + \frac{\lambda}{2\epsilon K} u^2 - B) d\boldsymbol{x},$$

By taking the time derivative of the new variable $U$, we then reformulate the system (2.5)-(2.6) as the following equivalent PDE system,

$$(3.6) \qquad \tau(\phi)\phi_t = -Z(\phi) U - \frac{\lambda}{\epsilon} T p'(\phi),$$

$$(3.7) \qquad U_t = \frac{1}{2} Z(\phi) \phi_t,$$

$$(3.8) \qquad u_t = D \Delta u + K p'(\phi) \phi_t,$$

where

$$(3.9) \qquad Z(\phi) = \frac{-\nabla \cdot (\kappa^2(\nabla \phi)\nabla \phi + \kappa(\nabla \phi)|\nabla \phi|^2 \boldsymbol{H}(\phi)) + \frac{1}{\epsilon^2} f(\phi)}{\sqrt{\frac{1}{2}|\kappa(\nabla \phi)\nabla \phi|^2 + \frac{1}{4\epsilon^2} F(\phi) + B}}.$$

The initial conditions are given by

$$\text{(3.10)} \quad \begin{cases} \phi(t=0) = \phi_0, \ u(t=0) = u_0, \\ U(t=0) = \sqrt{\frac{1}{2}|\kappa(\nabla\phi_0)\nabla\phi_0|^2 + \frac{1}{4\epsilon^2}F(\phi_0) + B}. \end{cases}$$

The boundary conditions are still (2.9) that we alluded before.

The system (3.6)-(3.8) also follows an energy dissipative law in terms of $\phi$, $u$ and the new variable $U$. By taking the $L^2$ inner product of (3.6) with $\phi_t$, of (3.7) with $-2U$, of (3.8) with $\frac{\lambda}{\epsilon K}T$, performing integration by parts and summing all equalities up, we can obtain the energy dissipation law of the new system (3.6)-(3.8) as

$$\text{(3.11)} \quad \frac{d}{dt}E(\phi, U, u) = -\|\sqrt{\tau(\phi)}\phi_t\|^2 - \frac{\lambda D}{\epsilon K}\|\nabla u\|^2 \leq 0.$$

Next we will develop time marching algorithms for solving the transformed system (3.6)-(3.8). The proposed schemes should formally follow the new energy dissipation law (3.11) in the discrete sense, instead of the energy law for the originated system (2.10). It can be shown that the discrete transformed energy is the approximation to the original energy with the corresponding order, which is verified by the rigorous error estimates and numerical simulations in [37].

Let $\delta t > 0$ be a time step size and set $t^n = n\delta t$ for $0 \leq n \leq N$ with $T = N\delta t$. We also denote the $L^2$ inner product of any two spatial functions $f_1(\boldsymbol{x})$ and $f_2(\boldsymbol{x})$ by $(f_1(\boldsymbol{x}), f_2(\boldsymbol{x})) = \int_\Omega f_1(\boldsymbol{x})f_2(\boldsymbol{x})d\boldsymbol{x}$, and the $L^2$ norm of the function $f(\boldsymbol{x})$ by $\|f\|^2 = (f, f)$. Let $\psi^n$ denotes the numerical approximation to $\psi(\cdot, t)|_{t=t^n}$ for any function $\psi$. We also define the following Sobolev spaces $H^1_{per}(\Omega) = \{\phi \text{ is periodic}, \phi \in H^1(\Omega)\}$.

We construct a second-order time marching scheme for solving the system (3.6)-(3.8) based on the second-order backward differentiation formula (BDF2), as follows.

**Scheme 2.** *Assuming $\phi^n, U^n, u^n$ and $\phi^{n-1}, U^{n-1}, u^{n-1}$ are known, we update $\phi^{n+1}, U^{n+1}, u^{n+1}$ by solving the following linear coupled system:*

$$\text{(3.12)} \quad \tau(\phi^{\star,n+1})\frac{3\phi^{n+1} - 4\phi^n + \phi^{n-1}}{2\delta t} + \frac{S_1}{\epsilon^2}(\phi^{n+1} - 2\phi^n + \phi^{n-1}) - S_2\Delta(\phi^{n+1} - 2\phi^n + \phi^{n-1})$$
$$= -Z^{\star,n+1}U^{n+1} - \frac{\lambda}{\epsilon}p'(\phi^{\star,n+1})u^{n+1},$$

$$\text{(3.13)} \quad 3U^{n+1} - 4U^n + U^{n-1} = \frac{1}{2}Z^{\star,n+1}(3\phi^{n+1} - 4\phi^n + \phi^{n-1}),$$

$$\text{(3.14)} \quad \frac{3u^{n+1} - 4u^n + u^{n-1}}{2\delta t} - D\Delta u^{n+1} = Kp'(\phi^{\star,n+1})\frac{3\phi^{n+1} - 4\phi^n + \phi^{n-1}}{2\delta t},$$

*where $Z^{\star,n+1} = Z(\phi^{\star,n+1})$ and $S_i, i = 1, 2$ are two positive stabilizing parameters.*

**Remark 3.1.** *When $S_1 = S_2 = 0$, the above scheme is the IEQ type scheme which had been developed in [2, 7, 32–36, 38–40]. Note even though the IEQ method is formally unconditionally energy stable, the spatial oscillations caused by the anisotropic coefficient $\kappa(\nabla\phi)$ can still make the scheme blow up for large time steps, which are illustrated in Fig. 4.1 and 4.5.*

Scheme (3.12)-(3.14) is totally linear since we handle the nonlinear terms by compositions of implicit and explicit discretization at $t^{n+1}$. Note that the new variable $U$ will not bring up extra computational cost since we can rewrite (3.13) as follows

$$(3.15) \qquad U^{n+1} = \frac{1}{2} Z^{\star,n+1} \phi^{n+1} + A_1,$$

where $A_1 = \frac{4U^n - U^{n-1}}{3} - \frac{1}{2} Z^{\star,n+1} \frac{4\phi^n - \phi^{n-1}}{3}$, and substitute (3.15) to (3.12), then the system (3.12)-(3.14) can be rewritten as

$$(3.16) \qquad Q(\phi^{n+1}) + \frac{\lambda}{\epsilon} p'(\phi^{\star,n+1}) u^{n+1} = f_1,$$

$$(3.17) \qquad \frac{\lambda}{\epsilon K} u^{n+1} - \frac{2\delta t}{3} \frac{\lambda}{\epsilon K} D \Delta u^{n+1} - \frac{\lambda}{\epsilon} p'(\phi^{\star,n+1}) \phi^{n+1} = f_2,$$

where

$$(3.18) \qquad Q(\phi^{n+1}) = \frac{3}{2\delta t} \tau(\phi^{\star,n+1}) \phi^{n+1} + \frac{S_1}{\epsilon^2} \phi^{n+1} - S_2 \Delta \phi^{n+1} + \frac{1}{2} Z^{\star,n+1} Z^{\star,n+1} \phi^{n+1},$$

and $f_1, f_2$ include only terms from previous time steps. In practice, we solve (3.16)-(3.17) directly to obtain $\phi^{n+1}$ and $u^{n+1}$, and then update the new variable $U^{n+1}$ from (3.15).

Furthermore, we notice $(Q(\phi), \psi) = (\phi, Q(\psi))$ when $\phi, \psi$ satisfy the boundary conditions in (2.9), that means the linear operator $Q(\phi)$ is self-adjoint. Moreover, $(Q(\phi), \phi) \geq 0$, where "=" is valid if and only if $\phi \equiv 0$, that means the linear operator $Q(\phi)$ is actually symmetric positive definite.

We now show the well-posedness of the weak form of the above linear system (3.16)-(3.17). In the following arguments, we will only consider the periodic boundary condition for convenience. For the case of no-flux boundary conditions, the proof is similar.

The weak form of (3.16)-(3.17) can be written as the following system with the unknowns $(\phi, u) \in (H^1_{per}, H^1_{per})(\Omega)$,

$$(3.19) \qquad \frac{3}{2\delta t}(\tau(\phi^{\star,n+1})\phi, \psi) + \frac{S_1}{\epsilon^2}(\phi, \psi) + S_2(\nabla \phi, \nabla \psi) + \frac{1}{2}(Z^{\star,n+1}\phi, Z^{\star,n+1}\psi)$$
$$+ \frac{\lambda}{\epsilon}(p'(\phi^{\star,n+1})u, \psi) = (f_1, \psi),$$

$$(3.20) \qquad \frac{\lambda}{\epsilon K}(u, v) + \frac{2\delta t}{3} \frac{\lambda D}{\epsilon K}(\nabla u, \nabla v) - \frac{\lambda}{\epsilon}(p'(\phi^{\star,n+1})\phi, v) = (f_2, v),$$

for any $(\psi, v) \in (H^1_{per}, H^1_{per})(\Omega)$.

We denote the above bilinear system (3.19)-(3.20) as

$$(3.21) \qquad (\boldsymbol{A}(\boldsymbol{X}), \boldsymbol{Y}) = (\boldsymbol{B}, \boldsymbol{Y}),$$

where $\boldsymbol{X} = (\phi, u)^T, \boldsymbol{Y} = (\psi, v)^T$ and $\boldsymbol{X}, \boldsymbol{Y} \in (H^1_{per}, H^1_{per})(\Omega)$.

**Theorem 3.1.** *The bilinear system* (3.21) *admits a unique solution* $(\phi, u) \in (H^1_{per}, H^1_{per})(\Omega)$.

*Proof.* (i) For any $\boldsymbol{X} = (\phi, u)^T$ and $\boldsymbol{Y} = (\psi, v)^T$ with $\boldsymbol{X}, \boldsymbol{Y} \in (H^1_{per}, H^1_{per})(\Omega)$, we have

$$(3.22) \qquad (\boldsymbol{A}(\boldsymbol{X}),) \boldsymbol{Y} \leq C_1(\|\phi\|_{H^1} + \|u\|_{H^1})(\|\psi\|_{H^1} + \|v\|_{H^1}),$$

where $C_1$ is a constant depending on $\delta t$, $S_1$, $S_2$, $\lambda$, $\epsilon$, $D$, $K$, $\|Z^{\star,n+1}\|_\infty$, and $\|p'(\phi^{\star,n+1})\|_\infty$. Therefore, the bilinear form $(\boldsymbol{A}(\boldsymbol{X}), \boldsymbol{Y})$ is bounded.

(ii) It is easy to derive that

$$
\begin{aligned}
(\boldsymbol{A}(\boldsymbol{X}), \boldsymbol{X}) &= \frac{3}{2\delta t}\|\sqrt{\tau(\phi^{\star,n+1})}\phi\|^2 + \frac{S_1}{\epsilon^2}\|\phi\|^2 + S_2\|\nabla\phi\|^2 \\
&\quad + \frac{1}{2}\|Z^{\star,n+1}\phi\|^2 + \frac{\lambda}{\epsilon K}\|u\|^2 + \frac{2\delta t \lambda D}{3\epsilon K}\|\nabla u\|^2 \\
&\geq C_2(\|\phi\|_{H^1}^2 + \|u\|_{H^1}^2),
\end{aligned}
\tag{3.23}
$$

where $C_2$ is a constant depending on $\delta t$, $S_1$, $S_2$, $\epsilon$, $\lambda$, $D$, $K$. Thus the bilinear form $(\boldsymbol{A}(\boldsymbol{X}), \boldsymbol{Y})$ is coercive.

Then from the Lax-Milgram theorem, we conclude the linear system (3.21) admits a unique solution $(\phi, u) \in (H^1_{per}, H^1_{per})(\Omega)$. $\square$

**Remark 3.2.** *If $S_2 = 0$, the coercivity of the bilinear form (3.19)-(3.20) in $H^1(\Omega)$ is not valid and we can only show the system admits a unique solution $\phi$ in $L^2(\Omega)$. Therefore, the $H^1$ bound of the numerical solution cannot be justified theoretically, that might bring up some substantial challenges for error estimates. On the contrary, when $S_2 \neq 0$, it is expected that optimal error estimates can be obtained without essential difficulties since the $H^1$ bound for $\phi$ is satisfied naturally. We will implement the subsequent error analysis in the future work by following the same lines as the analytical work for isotropic Allen-Cahn/Cahn-Hilliard models in [37].*

Now we prove the scheme (3.12)-(3.14) is unconditionally energy stable as follows.

**Theorem 3.2.** *The scheme (3.12)-(3.14) is unconditionally energy stable which satisfies the following discrete energy dissipation law,*

$$
\frac{1}{\delta t}(E^{n+1}_{bdf2} - E^n_{bdf2}) \leq -\|\sqrt{\tau(\phi^{\star,n+1})}\frac{3\phi^{n+1} - 4\phi^n + \phi^{n-1}}{2\delta t}\|^2 - \frac{\lambda D}{\epsilon K}\|\nabla u^{n+1}\|^2 \leq 0,
\tag{3.24}
$$

*where*

$$
\begin{aligned}
E^{n+1}_{bdf2} &= \frac{\|U^{n+1}\|^2 + \|2U^{n+1} - U^n\|^2}{2} + \frac{\lambda}{2\epsilon K}\Big(\frac{\|u^{n+1}\|^2 + \|2u^{n+1} - u^n\|^2}{2}\Big) \\
&\quad + \frac{S_1}{\epsilon^2}\frac{\|\phi^{n+1} - \phi^n\|^2}{2} + S_2\frac{\|\nabla\phi^{n+1} - \nabla\phi^n\|^2}{2}.
\end{aligned}
\tag{3.25}
$$

*Proof.* By taking the $L^2$ inner product of (3.12) with $\frac{1}{2}(3\phi^{n+1} - 4\phi^n + \phi^{n-1})$, and using the following identity,

$$
\frac{1}{2}(3a - 4b + c)(a - 2b + c) = \frac{1}{2}(a-b)^2 - \frac{1}{2}(b-c)^2 + (a - 2b + c)^2,
\tag{3.26}
$$

we obtain

$$
\begin{aligned}
\delta t \| \sqrt{\tau(\phi^{\star,n+1})} \frac{3\phi^{n+1} - 4\phi^n + \phi^{n-1}}{2\delta t} \|^2 & \\
+ \frac{S_1}{\epsilon^2}(\frac{1}{2}\|\phi^{n+1} - \phi^n\|^2 - \frac{1}{2}\|\phi^n - \phi^{n-1}\|^2 + \|\phi^{n+1} - 2\phi^n + \phi^{n-1}\|^2) & \\
+ S_2(\frac{1}{2}\|\nabla(\phi^{n+1} - \phi^n)\|^2 - \frac{1}{2}\|\nabla(\phi^n - \phi^{n-1})\|^2 + \|\nabla(\phi^{n+1} - 2\phi^n + \phi^{n-1})\|^2) & \\
= -\frac{1}{2}(Z^{\star,n+1}U^{n+1}, 3\phi^{n+1} - 4\phi^n + \phi^{n-1}) & \\
- \frac{\lambda}{2\epsilon}(p'(\phi^{\star,n+1})u^{n+1}, 3\phi^{n+1} - 4\phi^n + \phi^{n-1}). &
\end{aligned}
\tag{3.27}
$$

By taking the $L^2$ inner product of (3.13) with $U^{n+1}$ and using the following identity

$$
a(3a - 4b + c) = \frac{1}{2}a^2 + \frac{1}{2}(2a - b)^2 - \frac{1}{2}b^2 - \frac{1}{2}(2b - c)^2 + \frac{1}{2}(a - 2b + c)^2,
\tag{3.28}
$$

we obtain

$$
\begin{aligned}
(\frac{1}{2}\|U^{n+1}\|^2 + \frac{1}{2}\|2U^{n+1} - U^n\|^2) - (\frac{1}{2}\|U^n\|^2 + \frac{1}{2}\|2U^n - U^{n-1}\|^2) & \\
+ \frac{1}{2}\|U^{n+1} - 2U^n + U^{n-1}\|^2 = \frac{1}{2}(Z^{\star,n+1}(3\phi^{n+1} - 4\phi^n + \phi^{n-1}), U^{n+1}). &
\end{aligned}
\tag{3.29}
$$

By taking the $L^2$ inner product of (3.14) with $\delta t \frac{\lambda}{\epsilon K} u^{n+1}$, we obtain

$$
\begin{aligned}
\frac{\lambda}{2\epsilon K}(\frac{1}{2}\|u^{n+1}\|^2 + \frac{1}{2}\|2u^{n+1} - u^n\|^2) - \frac{\lambda}{2\epsilon K}(\frac{1}{2}\|u^n\|^2 + \frac{1}{2}\|2u^n - u^{n-1}\|^2) & \\
+ \frac{\lambda}{4\epsilon K}\|u^{n+1} - 2u^n + u^{n-1}\|^2 + \frac{\lambda D}{\epsilon K}\delta t \|\nabla u^{n+1}\|^2 & \\
= \frac{\lambda}{2\epsilon}(p'(\phi^{\star,n+1})(3\phi^{n+1} - 4\phi^n + \phi^{n-1}), u^{n+1}). &
\end{aligned}
\tag{3.30}
$$

By combing (3.27), (3.29) and (3.30), we obtain

$$
\begin{aligned}
(\frac{1}{2}\|U^{n+1}\|^2 + \frac{1}{2}\|2U^{n+1} - U^n\|^2) - (\frac{1}{2}\|U^n\|^2 + \frac{1}{2}\|2U^n - U^{n-1}\|^2) & \\
+ \frac{\lambda}{2\epsilon K}(\frac{1}{2}\|u^{n+1}\|^2 + \frac{1}{2}\|2u^{n+1} - u^n\|^2) - \frac{\lambda}{2\epsilon K}(\frac{1}{2}\|u^n\|^2 + \frac{1}{2}\|2u^n - u^{n-1}\|^2) & \\
+ \frac{S_1}{\epsilon^2}\frac{\|\phi^{n+1} - \phi^n\|^2}{2} - \frac{S_1}{\epsilon^2}\frac{\|\phi^n - \phi^{n-1}\|^2}{2} & \\
+ S_2\frac{\|\nabla(\phi^{n+1} - \phi^n)\|^2}{2} - S_2\frac{\|\nabla(\phi^n - \phi^{n-1})\|^2}{2}) & \\
+ \frac{1}{2}\|U^{n+1} - 2U^n + U^{n-1}\|^2 + \frac{\lambda}{4\epsilon K}\|u^{n+1} - 2u^n + u^{n-1}\|^2 & \\
+ \frac{S_1}{\epsilon^2}\|\phi^{n+1} - 2\phi^n + \phi^{n-1}\|^2 + S_2\|\nabla(\phi^{n+1} - 2\phi^n + \phi^{n-1})\|^2 & \\
= -\delta t \|\sqrt{\tau(\phi^{\star,n+1})}\frac{3\phi^{n+1} - 4\phi^n + \phi^{n-1}}{2\delta t}\|^2 - \frac{\lambda D}{\epsilon K}\delta t\|\nabla u^{n+1}\|^2. &
\end{aligned}
\tag{3.31}
$$

Finally, we obtain the desired result after dropping some positive terms. □

**Remark 3.1.** *Heuristically, $\frac{1}{\delta t}(E_{bdf2}^{n+1} - E_{bdf2}^{n})$ is a second-order approximation of $\frac{d}{dt}E(\phi, U, u)$ at $t = t^{n+1}$. For any smooth variable $\psi$ with time, we have*

$$\frac{\|\psi^{n+1}\|^2 - \|2\psi^{n+1} - \psi^n\|^2}{2\delta t} - \frac{\|\psi^n\|^2 - \|2\psi^n - \psi^{n-1}\|^2}{2\delta t}$$
$$\cong \frac{\|\psi^{n+2}\|^2 - \|\psi^n\|^2}{2\delta t} + O(\delta t^2) \cong \frac{d}{dt}\|\psi(t^{n+1})\|^2 + O(\delta t^2),$$

*and*

(3.32) $$\frac{\|\psi^{n+1} - \psi^n\|^2 - \|\psi^n - \psi^{n-1}\|^2}{2\delta t} \cong O(\delta t^2).$$

**Remark 3.2.** *It is also straightforward to develop the second-order Crank-Nicolson type scheme where the linear stabilizers terms still form like $\psi^{n+1} - 2\psi^n + \psi^{n-1}$. We leave the details to the interested readers since the proof of energy stability is quite similar to Theorem 3.2. In addition, although we consider only time discrete schemes in this study, the results can be carried over to any consistent finite-dimensional Galerkin approximations in the space since the proofs are all based on a variational formulation with all test functions in the same space as the space of the trial functions.*

## 4. Numerical simulations

In this section, we present various numerical examples to validate the proposed schemes and demonstrate their accuracy, energy stability and efficiency.

### 4.1. Accuracy test.
We first implement a numerical example with fourfold anisotropy (2.3) in 2D space to test the convergence rates of the two proposed schemes. We choose the periodic boundary conditions and set the computational domain as $\Omega = [0, h_1] \times [0, h_2]$. We use the Fourier-spectral method to discretize the space, where $129 \times 129$ Fourier modes are used. The model parameters are set as follows,

(4.1) $\quad h_1 = h_2 = 2\pi, \tau = 100, \epsilon = 0.06, \epsilon_4 = 0.05, D = \lambda = K = 1, S_1 = S_2 = 4.$

We assume the following two functions

(4.2) $\quad \phi(x, y, t) = \sin(x)\cos(y)\cos(t), \ u(x, y, t) = \cos(x)\sin(y)\cos(t)$

to be the exact solutions, and impose some suitable force fields such that the given solutions can satisfy the system (2.5)-(2.6). To see how the stability/accuracy is affected by the stabilizing terms, we compare the numerical results computed by four schemes, i.e., the stabilized linear scheme ((3.1)-(3.2), denoted by SLS) and its counterpart: the linear explicit scheme ((3.1)-(3.2) with $S_1 = S_2 = 0$, denoted by LS); the stabilized IEQ scheme ((3.12)-(3.14), denoted by SIEQ) and its counterpart: the IEQ scheme ((3.12)-(3.14) with $S_1 = S_2 = 0$, denoted by IEQ).

In Fig. 4.1, we plot the $L^2$ errors of the variables $\phi$ and $u$ between the numerically simulated solution and the exact solution at $T = 1$ with different time step sizes. Some remarkable features observed from Fig. 4.1 are listed as follows.

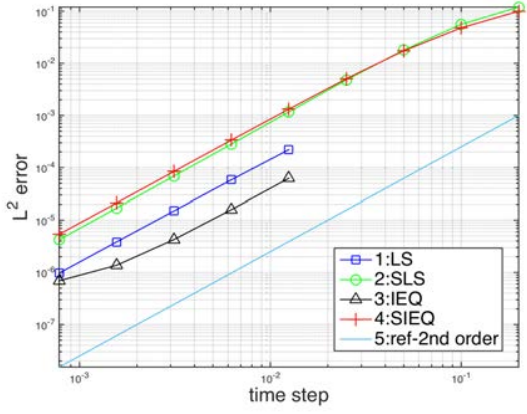 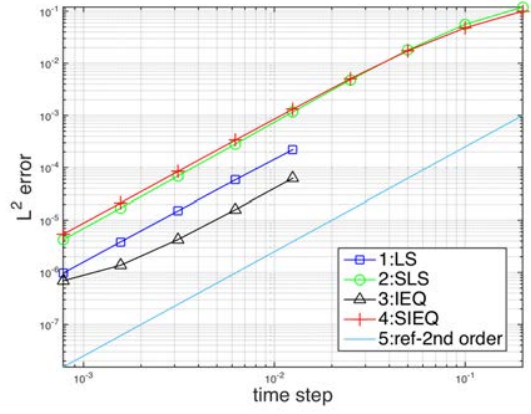

(a) The $L^2$ errors for $\phi$ with various time steps.  (b) The $L^2$ errors for $u$ with various time steps.

FIGURE 4.1. The $L^2$ numerical errors for the phase variable $\phi$ and temperature $u$ at $T = 1$, that are computed using the schemes LS, SLS, IEQ and SIEQ and various temporal resolutions with the initial conditions of (4.2).

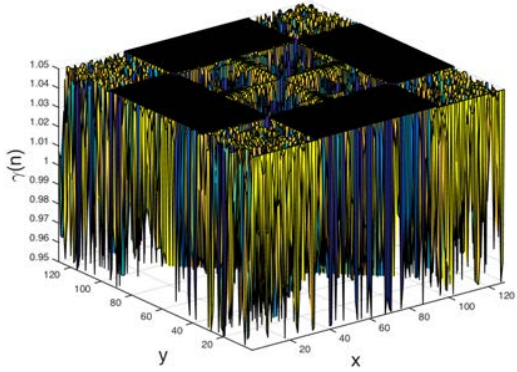 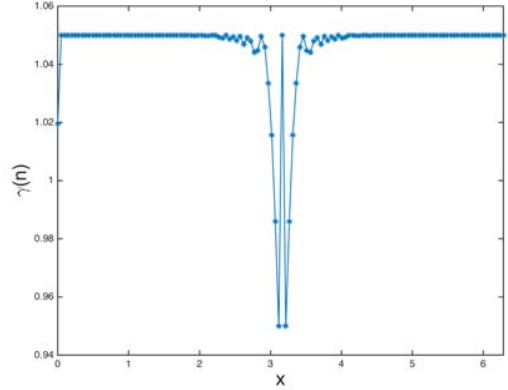

FIGURE 4.2. The profile of $\gamma(\bm{n}_0)$ with $\epsilon_4 = 0.25$ and the initial condition (4.3). The left subfigure is the 2D surface plots of $\gamma(\bm{n}_0)$, and the right subfigure is the 1D cross-section of $\gamma(\bm{n}_0)|_{(\cdot, y=\pi)}$.

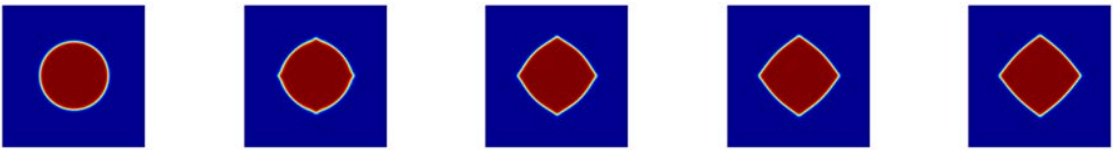

FIGURE 4.3. The 2D dynamical evolution of the phase variable $\phi$ by using the initial condition (4.3) and the time step $\delta t = 1e-2$. Snapshots of the numerical approximation are taken at $t = 0, 5, 10, 15,$ and $20$.

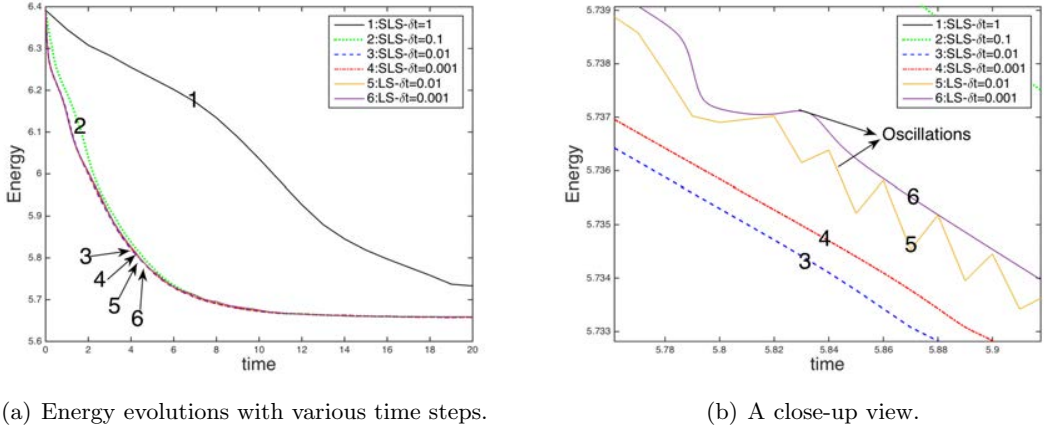

(a) Energy evolutions with various time steps.

(b) A close-up view.

FIGURE 4.4. Time evolutions of the total free energy functional using the two schemes, LS and SLS, for four different time steps $\delta t = 1$, $1e-1$, $1e-2$, and $1e-3$ (the LS scheme blows up quickly for $\delta t = 1$ and $1e-1$, thus the two curves are missing). The left subfigure (a) is the energy profile for $t \in [0, 20]$, and the right subfigure (b) is a close-up view for $t \in [5.7, 5.95]$.

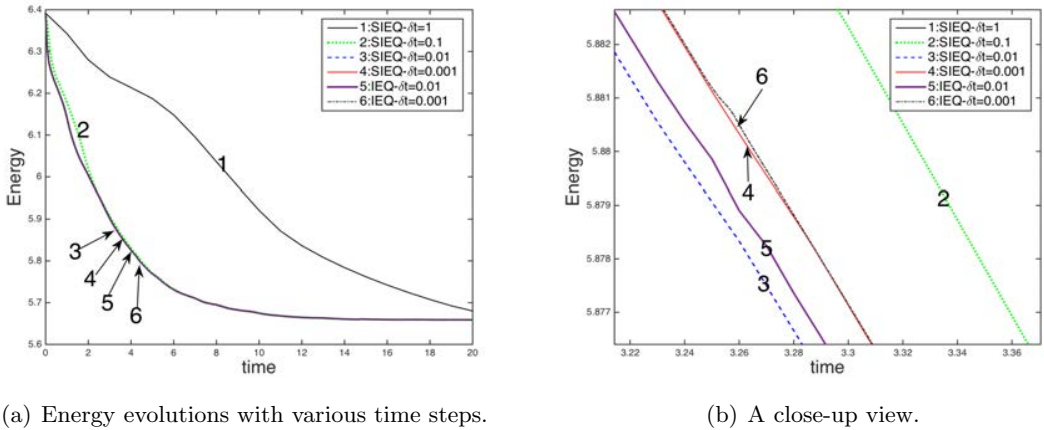

(a) Energy evolutions with various time steps.

(b) A close-up view.

FIGURE 4.5. Time evolutions of the total free energy functional using the two schemes, IEQ and SIEQ, for four different time steps $\delta t = 1$, $1e-1$, $1e-2$, and $1e-3$ (the IEQ scheme blows up quickly for $\delta t = 1$ and $1e-1$ thus the two energy curves are missing). The left subfigure (a) is the energy profile for $t \in [0, 20]$, and the right subfigure (b) is a close-up view for $t \in [3.21, 3.37]$.

- When $\delta t > 1.25e-2$, the schemes LS and IEQ blow up quickly, thus the error points are missing in Fig. 4.1. But their stabilized versions, SLS and SIEQ, still work well and present the second-order accuracy.

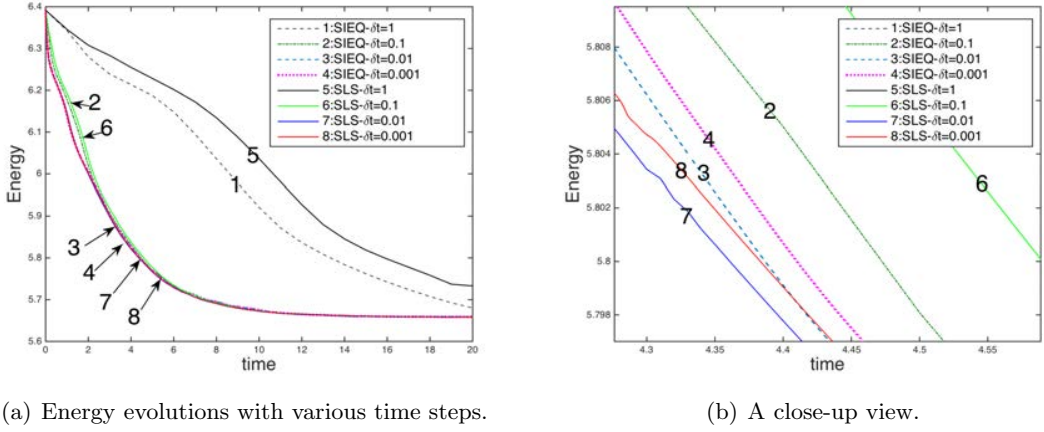

(a) Energy evolutions with various time steps.

(b) A close-up view.

FIGURE 4.6. Time evolutions of the total free energy functionals using the two stabilized schemes, SIEQ and SLS, for four different time steps $\delta t = 1$, 1e−1, 1e−2, and 1e−3. The left subfigure (a) is the energy profile for $t \in [0, 20]$, and the right subfigure (b) is a close-up view for $t \in [4.25, 4.6]$.

- When $\delta t \leq 1.25\text{e}{-2}$, all four schemes achieve almost perfect second-order accuracy in time. But obviously, the magnitude of errors computed by the SLS and SIEQ schemes are bigger than that computed by the LS and IEQ schemes. This phenomenon is reasonable since the added stabilizers actually increase the splitting errors indeed.

Therefore, in this example, when concerning the stability, the stabilized schemes (SLS and SIEQ) overwhelmingly defeat their non-stabilized counterparts (LS and IEQ). Meanwhile, in Fig. 4.1, we cannot tell obvious differences between the two stabilized schemes, SLS and SIEQ, from any point of view of the stability and/or accuracy.

4.2. **Evolution of a circle.** In this subsection, we consider the evolution of a circle driven by the strong fourfold anisotropy where the magnitude of $\epsilon_4$ is increased to be 0.25 and all other parameters are still from (4.1). The initial conditions for $\phi$ and $u$ read as

$$(4.3) \quad \begin{cases} \phi(x, y, t = 0) = \tanh(\dfrac{r_0 - \sqrt{(x - x_0)^2 + (y - y_0)^2}}{\epsilon_0}), \\ u(x, y, t = 0) = u_0, \end{cases}$$

where $(x_0, y_0, r_0, \epsilon_0, u_0) = (\pi, \pi, 1.5, 0.072, -0.55)$. We use the Fourier-spectral method to discretize the space, where $129 \times 129$ Fourier modes are used.

In Fig. 4.2, we present the 2D profile of $\gamma(\boldsymbol{n}_0)$ and 1D cross-section of $\gamma(\boldsymbol{n}_0)|_{(\cdot, y=\pi)}$, where one can observe high oscillations appear almost everywhere. We use the SIEQ scheme with the time step $\delta t = 1\text{e}{-2}$ to perform simulations. In Fig. 4.3, we show the dynamics how a circular shape interface with full orientations evolves to an anisotropic diamond shape with missing orientations at four corners. Snapshots of the phase field variable $\phi$ are taken at $t = 0, 5, 10, 15$, and $20$.

In next three figures, Fig. 4.4, 4.5, and 4.6, we compare the energy evolution curves computed by SLS and LS; SIEQ and IEQ; SLS and SIEQ, with various time steps, respectively. Some remarkable features are listed below.

- In Fig. 4.4, we compare the evolution of the free energy functionals using the schemes SLS and LS until the $T = 20$ with four different time steps $\delta t = 1$, 1e−1, 1e−2, and 1e−3. When $\delta t = 1$ and 1e−1, the LS scheme blows up quickly, thus two energy curves are missing. Moreover, even for small time steps $\delta t = $ 1e−2 and 1e−3, the discrete energies computed by LS still present some non-physical oscillations, shown in Fig. 4.4(b). But after we add the two stabilizers, these oscillations are removed and energies always decay monotonically, which means that the two stabilizers can suppress high-frequency oscillations efficiently.
- In Fig. 4.5, we further compare the evolution of the free energy functionals using the schemes SIEQ and IEQ until the $T = 20$ with four different time steps $\delta t = 1$, 1e−1, 1e−2, and 1e−3. For larger times steps $\delta t = 1$ and 1e−1, the IEQ scheme also blows up. But when we add two stabilizers in, all four energy curves generated by scheme SIEQ can decay monotonically, which numerically confirms that the SIEQ scheme is unconditionally energy stable for all time step sizes.
- In Fig. 4.6, we compare the energy evolution curves computed by the two stabilized schemes: SIEQ and SLS, with the four different time steps $\delta t = 1$, 1e−1, 1e−2, and 1e−3. We observe that, for the largest time step $\delta t = 1$, both of the energy curves computed by SLS and SIEQ deviate viewable away from others computed with smaller time steps. When using smaller time steps $\delta t = $ 1e−1, 1e−2, and 1e−3, the energy curves computed by SIEQ and SLS coincide very well. This means the adopted time step size for this example should not be larger than 1e−1, in order to get reasonably good accuracy.

4.3. **2D dendrite crystal growth with fourfold anisotropy.** In this subsection, we investigate how the anisotropic entropy coefficient can affect the shape of the dendritic crystal through the dynamical process in which a small crystal nucleus grows heterogeneously in 2D space. We set the initial condition as

(4.4)
$$\begin{cases} \phi(x, y, t = 0) = \tanh(\frac{r_0 - \sqrt{(x-x_0)^2 + (y-y_0)^2}}{\epsilon_0}); \\ u(x, y, t = 0) = \begin{cases} 0, & \phi > 0; \\ u_0, & \text{otherwise}, \end{cases} \end{cases}$$

where $(x_0, y_0, r_0, \epsilon_0, u_0) = (\pi, \pi, 0.02, 0.072, -0.55)$. The other parameters are set as follows,

(4.5)
$$\begin{aligned} &h_1 = h_2 = 2\pi, \quad \tau = 4.4\text{e}3, \quad \epsilon = 1.12\text{e}{-2}, \quad \epsilon_4 = 0.05, \\ &D = 2.25\text{e}{-4}, \quad K = 0.5, \quad \lambda = 380, \quad\quad S_1 = S_2 = 4. \end{aligned}$$

Note these parameters are the rescaled values obtained from the literatures [12, 14, 15, 40]. We use the Fourier-spectral method to discretize the space, where $513 \times 513$ Fourier modes are used.

We perform a series of simulations with the fourfold anisotropic entropy coefficient (2.3) by varying the latent heat parameter $K$. For all simulations, we use the scheme SIEQ (3.12)-(3.14) since it

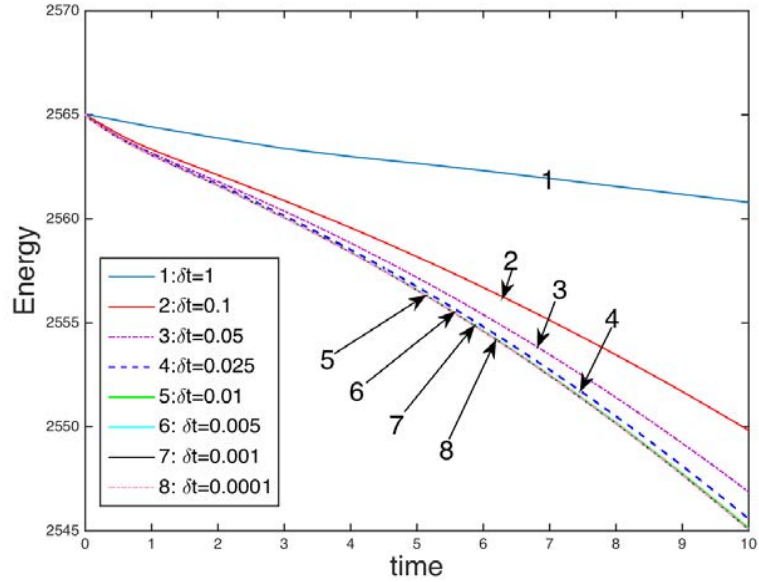

Figure 4.7. Time evolutions of the total free energy functional using the SIEQ scheme, for eight different time steps $\delta t = 1$, 0.1, 0.05, 0.025, 0.01, 0.005, 0.001, and 0.0001 with the initial condition (4.6) and order parameters (4.5). The energy curves show the decays for all time step sizes, that confirms that our algorithm is unconditionally stable for all tested time steps.

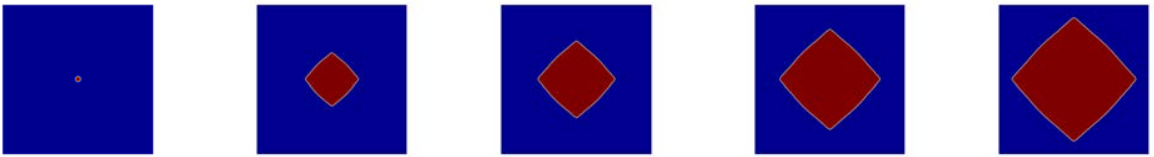

(a) The profiles of the phase field variable $\phi$.

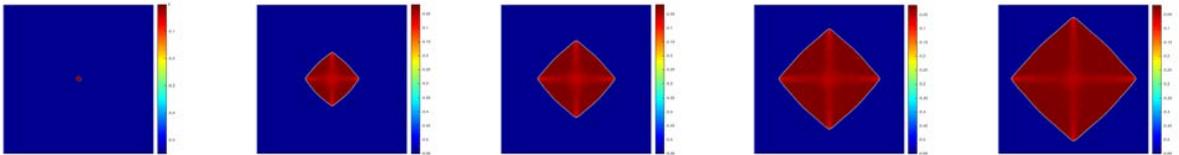

(b) The profiles of the temperature variable $u$.

Figure 4.8. The 2D dynamical evolutions of the dendritic crystal growth process with $K = 0.5$ and default parameters (4.5), computed by the scheme SIEQ and the time step $\delta t = 1e-2$. Snapshots of the numerical approximation are taken at $T = 0$, 40, 60, 80, and 100.

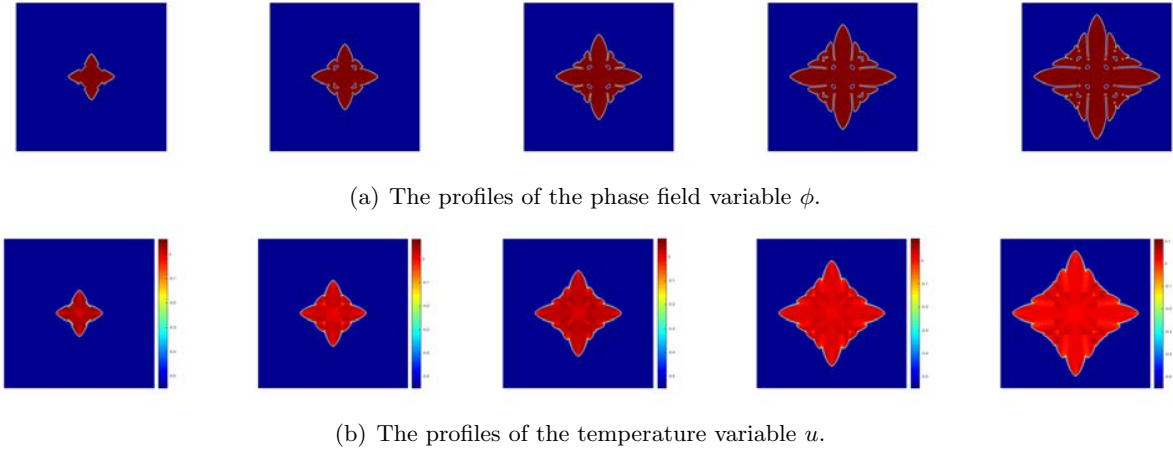

(a) The profiles of the phase field variable $\phi$.

(b) The profiles of the temperature variable $u$.

FIGURE 4.9. The 2D dynamical evolution of dendritic crystal growth process with $K = 0.6$ and order parameters (4.5), computed by the scheme SIEQ and the time step $\delta t = 1\mathrm{e}{-2}$. Snapshots of the numerical approximation are taken at $T = 40$, 60, 80, 100, and 120.

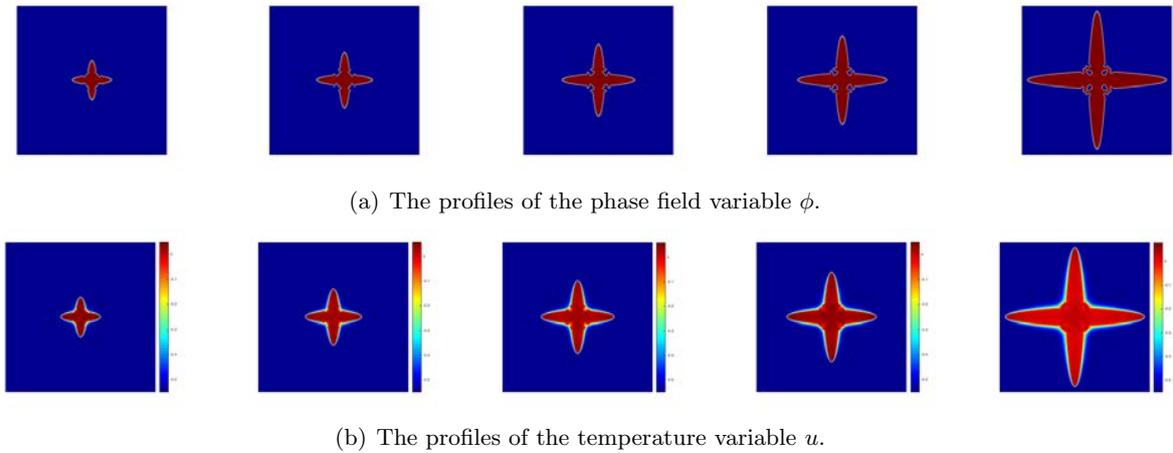

(a) The profiles of the phase field variable $\phi$.

(b) The profiles of the temperature variable $u$.

FIGURE 4.10. The 2D dynamical evolution of dendritic crystal growth process with $K = 0.7$ and order parameters (4.5), computed by the scheme SIEQ and the time step $\delta t = 1\mathrm{e}{-2}$. Snapshots of the numerical approximation are taken at $T = 40$, 60, 80, 100, and 160.

is provably energy stable. However, the stable scheme allows for any time step in computations only from the stability concern instead of accuracy. It is obvious that larger time step can lead to large computational errors. Therefore, to obtain good accuracy while consuming as low computational cost as possible, one needs to discover the rough range of the allowable maximum time step, which can be found through the energy evolution curve plots, shown in Fig. 4.7, where we compare the

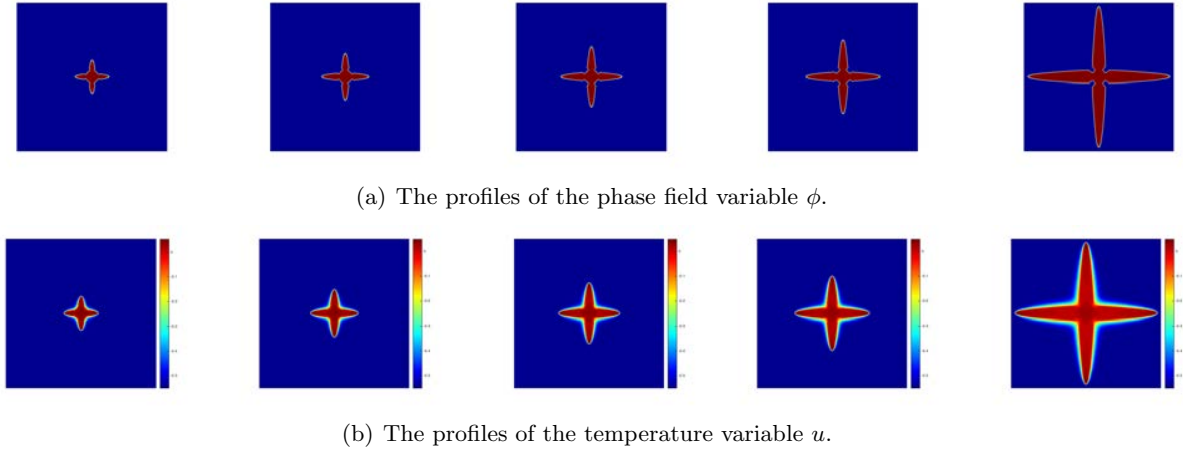

(a) The profiles of the phase field variable $\phi$.

(b) The profiles of the temperature variable $u$.

FIGURE 4.11. The 2D dynamical evolution of dendritic crystal growth process with $K = 0.8$ and default parameters (4.5), computed by the scheme SIEQ and the time step $\delta t = 1\text{e}-2$. Snapshots of the numerical approximation are taken at $T = 40$, 60, 80, 100, and 200.

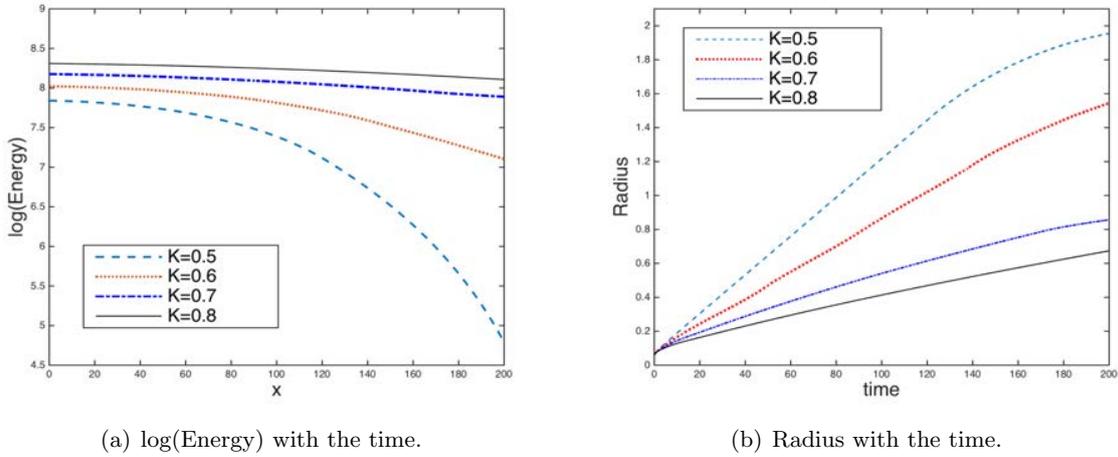

(a) log(Energy) with the time.

(b) Radius with the time.

FIGURE 4.12. (a) Time evolutions of the logarithm of the free energy functional for the 2D fourfold examples with $K = 0.5$, 0.6, 0.7, and 0.8. (b) The size of the crystals changing with time where the crystal size is measured by an equivalent radius of a circle with the same area.

time evolution of the free-energy functional for eight different time step sizes until $T = 10$. Fig. 4.7 delivers three messages as follows: (i) all eight tested energy curves show decays for all time step sizes, which verifies that the proposed SIEQ algorithm is unconditionally energy stable; (ii) the solution computed with smaller time steps of $\delta t = 0.01$, 0.005, 0.001, and 0.0001 can be treated

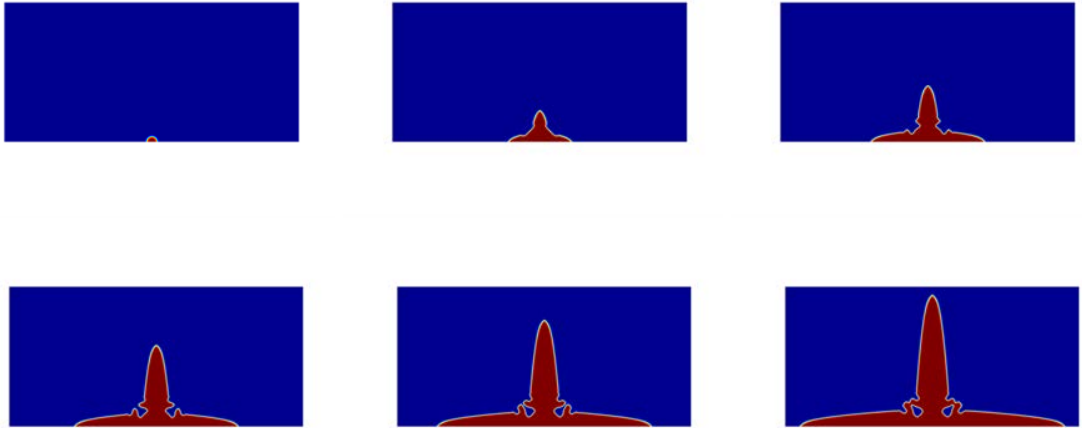

(a) The profiles of the phase field variable $\phi$.

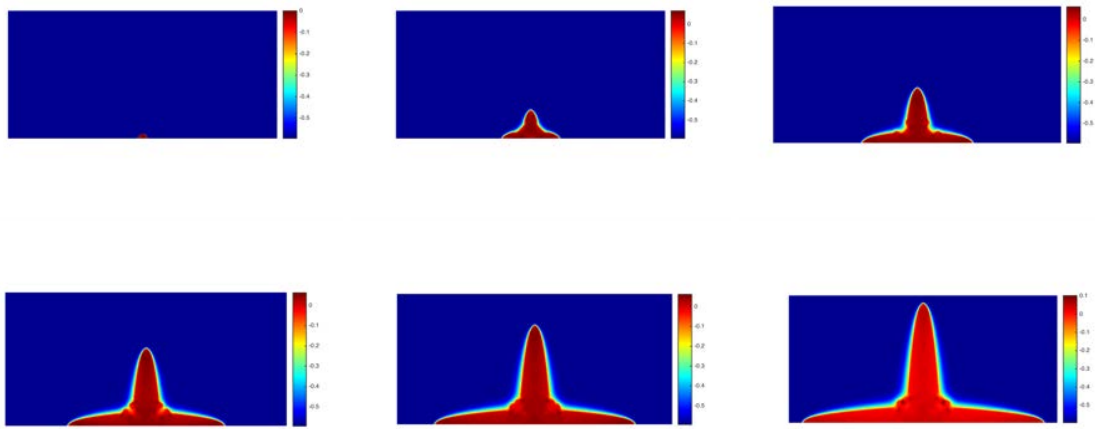

(b) The profiles of the temperature $u$.

FIGURE 4.13. The 2D dynamical evolution of dendritic crystal growth process, where the initial nucleus is set at the bottom with the no-flux boundary conditions for $y-$direction, computed by using the scheme SIEQ, with the time step $\delta t = 1\mathrm{e}-2$ and default parameters (4.5). Snapshots of the numerical approximation are taken at $T = 0, 30, 60, 90, 120,$ and $150$.

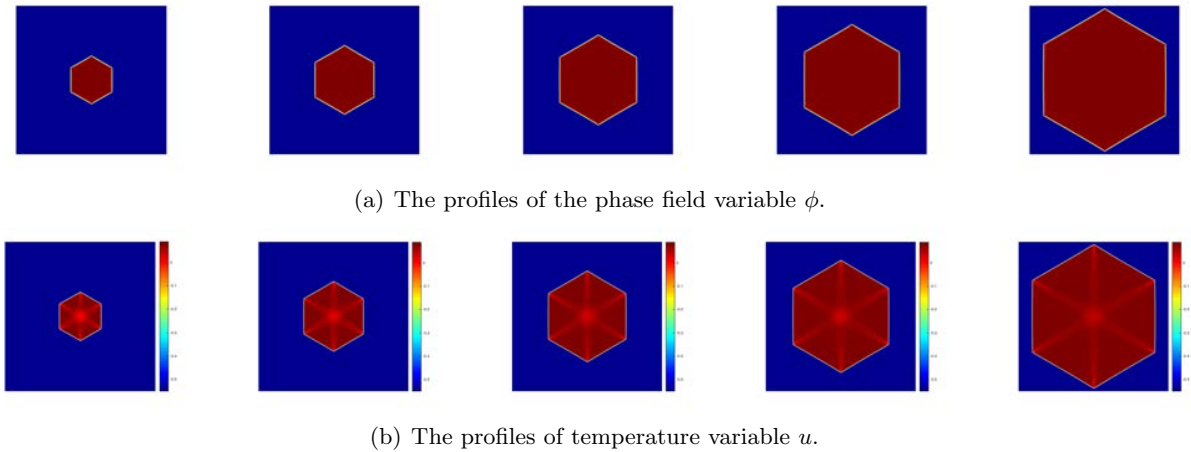

(a) The profiles of the phase field variable $\phi$.

(b) The profiles of temperature variable $u$.

FIGURE 4.14. The 2D dynamical evolution of dendritic crystal growth process with sixfold anisotropy and $K = 0.6$, computed by using the scheme SIEQ, the time step $\delta t = 1\mathrm{e}{-2}$, and default parameters (4.5). Snapshots of the numerical approximation are taken at $T = 40, 60, 80, 100$, and $130$.

as benchmark solutions since the four energy curves coincide very well; (iii) when $\delta t \geq 0.025$, the energy curves deviate viewable away from the benchmark solution, which means the time step size has to be smaller than 0.025 at least, in order to get reasonable accuracy. Therefore, in all following 2D simulations, we choose $\delta t = 0.01$.

In Fig. 4.8(a), we show snapshots of the phase variable at various times with $K = 0.5$. The tiny circle at the initial moment $t = 0$, shown in the first subfigure in Fig. 4.8(a), works as a crystal nucleus to grow with the time. Due to the anisotropic effects, we then observe the growth of the crystalline phase that finally becomes an anisotropic shape with missing orientations at four corners. In Fig. 4.8(b), we show the profiles of the temperature field $u$, that agrees well with the phase field variable $\phi$ due to the latent heat coupling terms.

Next, we adjust the magnitude of the latent heat parameter $K$ such that $K = 0.6$, $0.7$, and $0.8$ in Fig. 4.9, 4.10, and 4.11, respectively. To compare with the results of $K = 0.5$ in Fig. 4.8, all other parameters are kept to be the same as the default values in (4.5), and the same time step $\delta t = 0.01$ is used for better accuracy. When $K = 0.6$, in Fig. 4.9, we observe that the crystal nucleus initially grows to a star-shape where four sharper tips appear ($T = 40$), and finally grows like a snowflake pattern ($T = 120$) but with four branches. After we further increase $K$ to be $0.7$ and $0.8$ in Fig. 4.10 and Fig. 4.11, respectively, we then observe sharper tips and thinner branches formed. All these numerical results demonstrate similar features to those obtained in [13–15, 40]. In Fig. 4.12, for $K = 0.5, 0.6, 0.7$, and $0.8$, we summarize the evolutions of the logarithm of the total free energy which monotonically decay and the radius of the crystal that is measured by an equivalent radius of a circle with the same area.

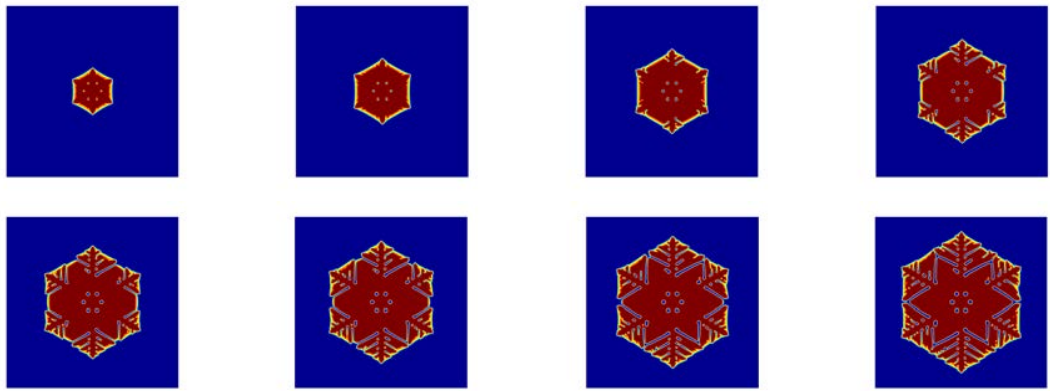

(a) The profile of the phase field variable $\phi$.

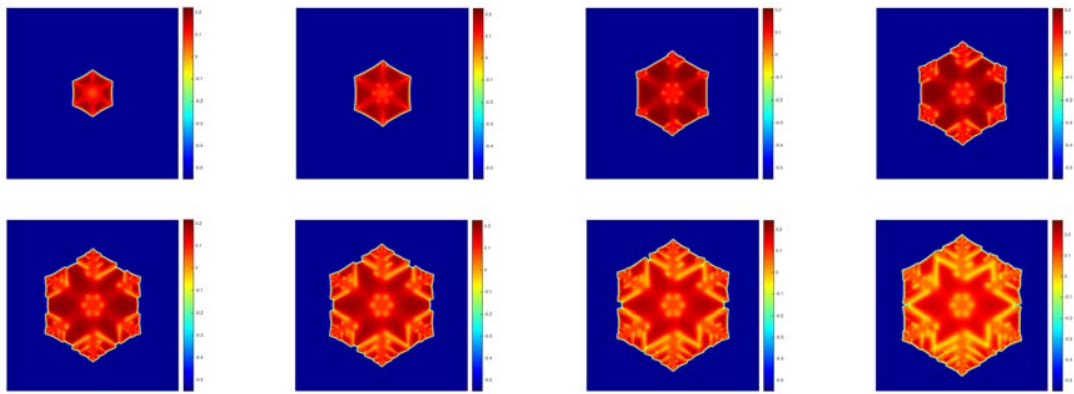

(b) The profiles of temperature variable $u$.

FIGURE 4.15. The 2D dynamical evolution of dendritic crystal growth process with sixfold anisotropy and $K = 0.7$, computed by using the scheme SIEQ, the time step $\delta t = 1\mathrm{e}{-2}$ and default parameters (4.5). Snapshots of the numerical approximation are taken at $T = 40, 60, 80, 100, 110, 120, 130$, and $140$.

We further perform a numerical test with the no-flux boundary condition ((2.9)-(ii)). In Fig. 4.13, we set the $y$-direction to be no-flux and still keep the $x$-direction to be periodic. We set $K = 0.75$, $\delta t = 0.01$ and the computed domain is $[0, h_1] \times [0, h_2]$ where $h_1 = 2\pi$ and $h_2 = 3$. The initial nucleus is set at the bottom (i.e., $y_0 = 0$ in (4.6)), shown in the first subfigure of Fig. 4.13(a). In Fig. 4.13(a), the crystal nucleus grows to the half of the snowflake pattern with four branches, which are consistent to the profiles shown in Fig. 4.10 with periodic boundary conditions. The dynamical evolutions of the temperature $u$ are shown in Fig. 4.13(b).

4.4. **2D dendrite crystal growth with with sixfold anisotropy.** In this subsection, we consider the alternative anisotropy coefficient (2.2) with $m = 6$ to investigate how the sixfold anisotropy can affect the shape of the dendritic crystal. To compare with the fourfold case, we use the same

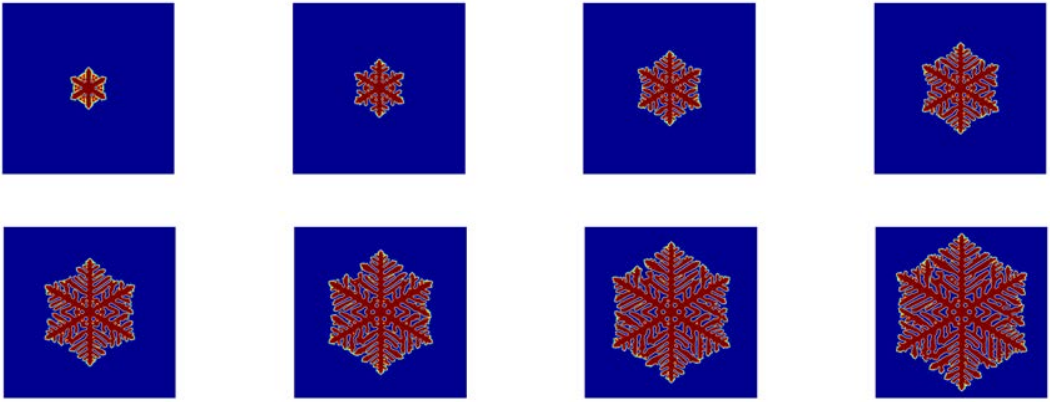

(a) The profiles of the phase field variable $\phi$.

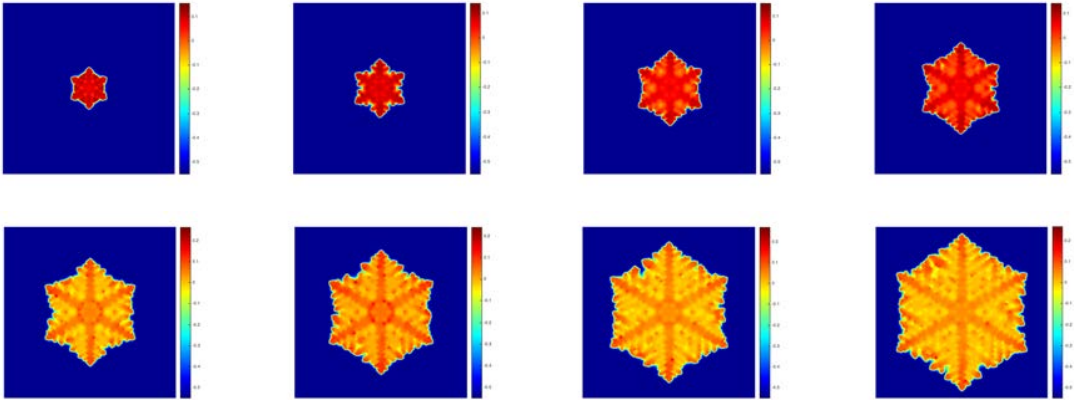

(b) The profiles of temperature variable $u$.

FIGURE 4.16. The 2D dynamical evolutions of dendritic crystal growth process with sixfold anisotropy and $K = 0.75$, computed by using the scheme SIEQ, the time step $\delta t = 1e-2$ and default parameters (4.5). Snapshots of the numerical approximation are taken at $T = 40, 60, 80, 100, 120, 140, 160,$ and $180$.

initial condition (4.6) and order parameters (4.5). The scheme SIEQ with $513 \times 513$ Fourier modes are used, and the time step $\delta t$ is set to be 0.01 for better accuracy.

We set the latent heat parameter $K$ to be 0.6, 0.7, 0.75, and 0.8 in Fig. 4.14, 4.15, 4.16, and 4.17, respectively. In Fig. 4.14 with $K = 0.6$, the initial circular nucleus grows into the regular hexagon shape. When $K = 0.7$, in Fig. 4.15, although the nucleus finally forms the hexagon shape, plenty of subtle microstructures are formed due to the anisotropy in the heat transfer process. When $K = 0.75$, in Fig. 4.16, the nucleus initially grows up with six main branches. When time evolves, plenty of sub-branches form on each main branch and finally, a snowflake pattern is presented. Similar snowflake patterns had been reported in [15] using a slightly different model.

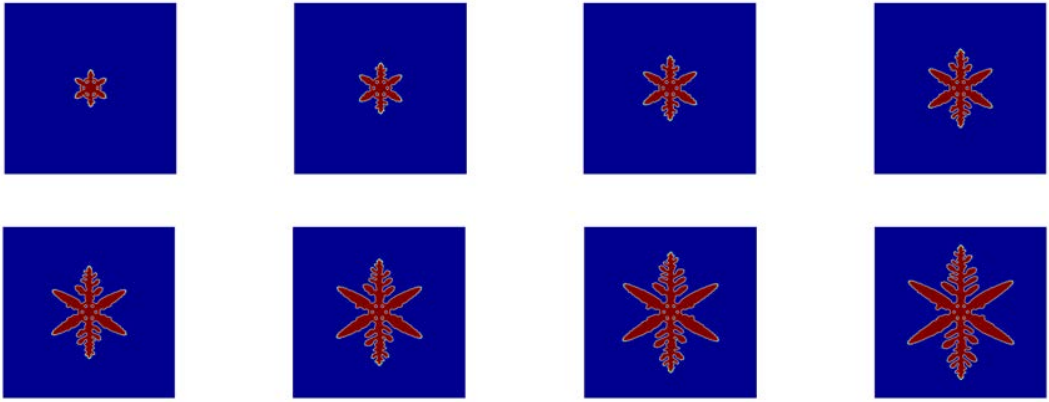

(a) The profile of the phase field variable $\phi$.

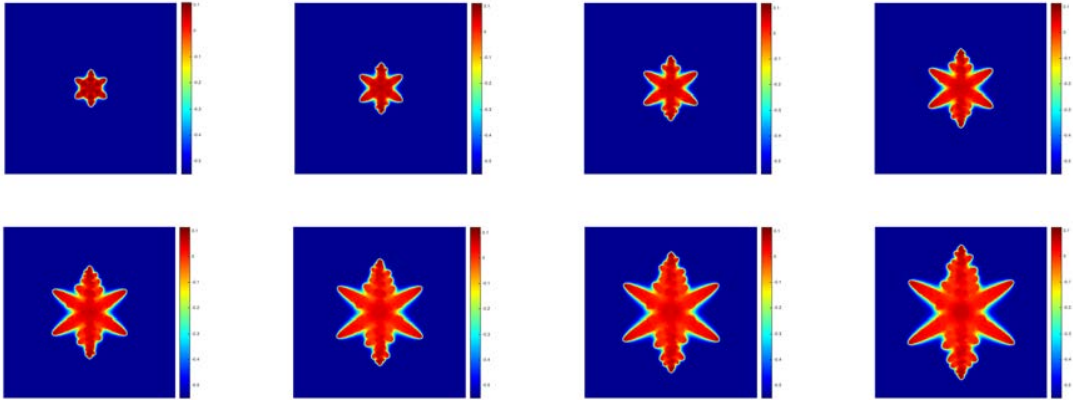

(b) The profiles of temperature variable $u$.

FIGURE 4.17. The 2D dynamical evolution of dendritic crystal growth process with sixfold anisotropy and $K = 0.8$, computed by using the scheme SIEQ, the time step $\delta t = 1\text{e}{-2}$ and default parameters (4.5). Snapshots of the numerical approximation are taken at $T = 40, 60, 80, 100, 120, 140, 160,$ and $180$.

When $K = 0.8$, in Fig. 4.17, the snowflake pattern but with fewer sub-branches are formed. In Fig. 4.18, we summarize the evolutions of the total free energy and crystal sizes for these four cases.

Finally, we set $m = 5$ for the fivefold anisotropy and $m = 7$ for the sevenfold anisotropy in Fig. 4.19(a) and Fig. 4.19(b), respectively. From the shapes of each $m$, we conclude that the number of dendrite branches is determined by $m$ exclusively. In the last subfigure of each figure, we show the corresponding temperature field $u$.

4.5. **3D dendrite crystal growth with fourfold anisotropy.** In this example, we investigate how the anisotropic coefficient can affect the shape of the dendritic crystal in 3D space. We set the

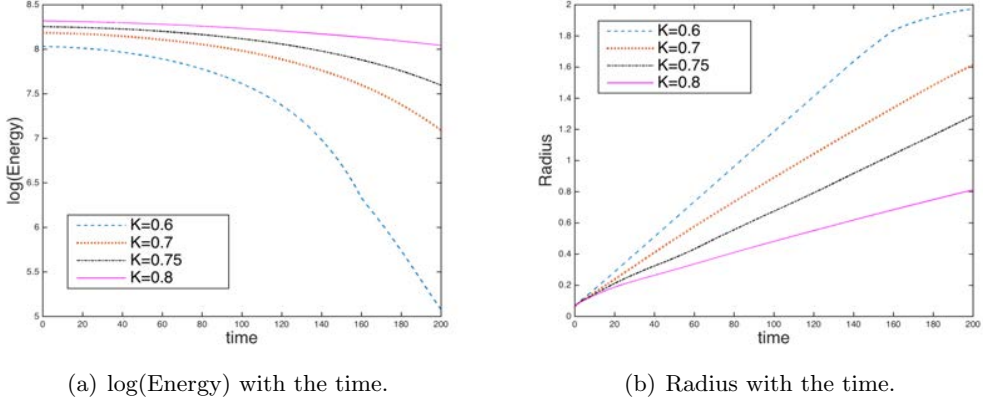

(a) log(Energy) with the time.  (b) Radius with the time.

FIGURE 4.18. (a) Time evolutions of the free energy functional for the sixfold examples with $K = 0.6$, $0.7$, $0.75$, and $0.8$. (b) The size of the crystals changing with time where the crystal size is measured by an equivalent radius of a circle with the same area.

initial condition as

(4.6)
$$\begin{cases} \phi(x,y,z,t=0) = \tanh(\dfrac{r_0 - \sqrt{(x-x_0)^2 + (y-y_0)^2 + (z-z_0)^2}}{\epsilon_0}); \\ u(x,y,z,t=0) = \begin{cases} 0, & \phi > 0; \\ u_0, & \text{otherwise}, \end{cases} \end{cases}$$

where $(x_0, y_0, z_0, r_0, \epsilon_0, u_0) = (\pi, \pi, \pi, 0.02, 0.072, -0.55)$. The other parameters are set as follows,

(4.7)
$$\begin{aligned} h_1 = h_2 = 2\pi, \quad \tau = 2.5\text{e}4, \quad \epsilon = 3\text{e}{-2}, \quad \epsilon_4 = 0.05, \\ D = 2\text{e}{-4}, \quad \lambda = 260, \quad S_1 = S_2 = 4. \end{aligned}$$

We use the Fourier-spectral method to discretize the space, where $129 \times 129 \times 129$ Fourier modes are used.

We use the scheme SIEQ to perform the 3D simulations with the fourfold anisotropic entropy coefficient (2.4), and use the time step $\delta t = 0.1$. We set $K = 0.5$ in Fig. 4.20, where we observe that the initial tiny 3D nucleus evolves to an anisotropic pyramid with missing orientations at six corners. In Fig. 4.21, 4.22, and 4.23 with $K = 1$, $1.5$, and $2$, respectively, we observe six branches form. Meanwhile, while $K$ increases, the formed tips and branches become sharper and thinner. This can be viewed more clearly in Fig. 4.24, where we show the 2D cross-section $\phi(\pi, \cdot, \cdot)$ for the computed solutions in Fig. 4.20-4.23 with $T = 900$, $1500$, $2100$, and $2700$, respectively. In Fig. 4.25, we summarize the evolutions of the total free energy and crystal sizes for the cases of $K = 0.5$, $1$, $1.5$, $2$, respectively. These 3D results are consistent to the 2D simulations with fourfold anisotropy shown in Section 4.3.

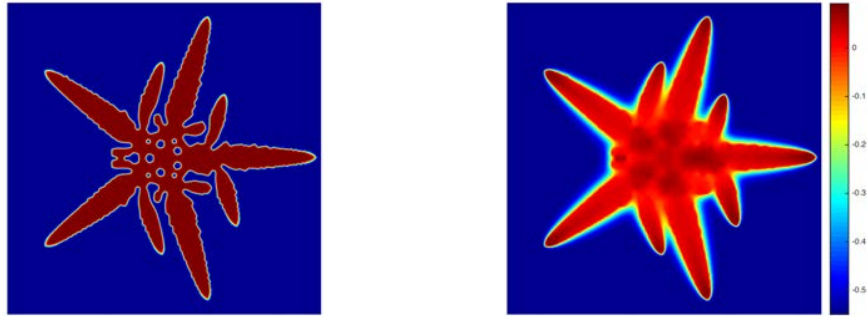

(a) The profiles of phase field $\phi$ (left), and the temperature $u$ (right) for fivefold anisotropy.

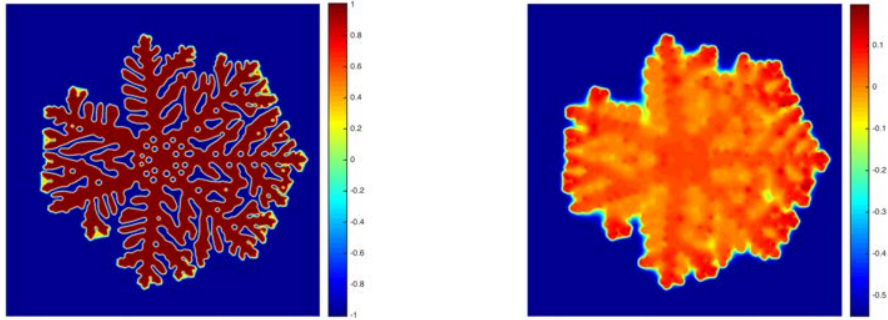

(b) The profiles of phase field $\phi$ (left), and the temperature $u$ (right) for sevenfold anisotropy.

FIGURE 4.19. The formed dendritic crystal with fivefold and sevenfold anisotropy, computed by using the scheme SIEQ, the time step $\delta t = 1e-2$ and default parameters (4.5). For each subfigure, the left one is the profile of the phase field variable $\phi$, and the right one is the temperature $u$ at the same moment. Snapshots of the numerical approximation are taken at $T = 170$ for (a), and $T = 200$ for (b).

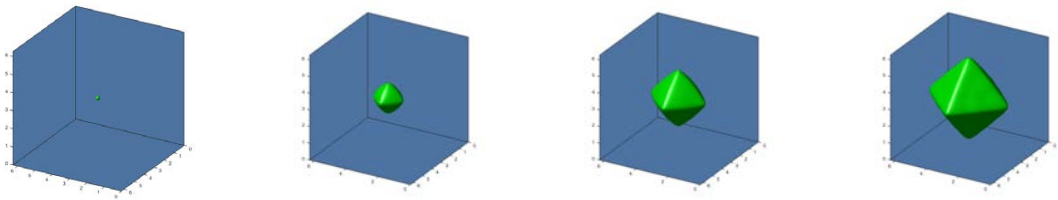

FIGURE 4.20. The dynamical evolution of dendritic crystal growth process in 3D with fourfold anisotropy and $K = 0.5$, computed by using the scheme SIEQ, the time step $\delta t = 1e-1$, and parameters (4.7). Snapshots of the numerical approximation are taken at $T = 0, 300, 600,$ and $900$.

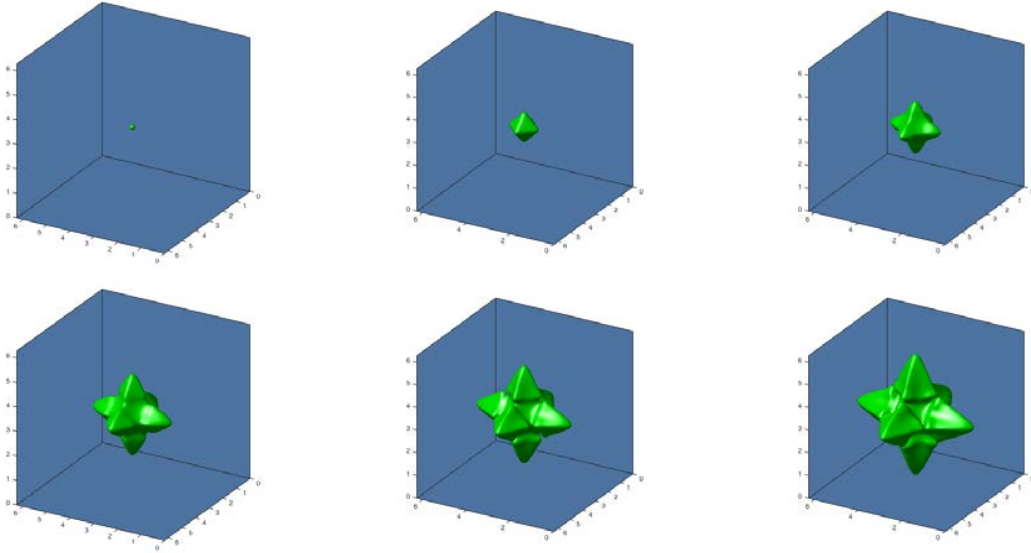

Figure 4.21. The dynamical evolution of dendritic crystal growth process in 3D with fourfold anisotropy and $K = 1$, computed by using the scheme SIEQ, the time step $\delta t = 1\mathrm{e}{-1}$, and parameters (4.7). Snapshots of the numerical approximation are taken at $T = 0$, 300, 600, 900, 1200, and 1500.

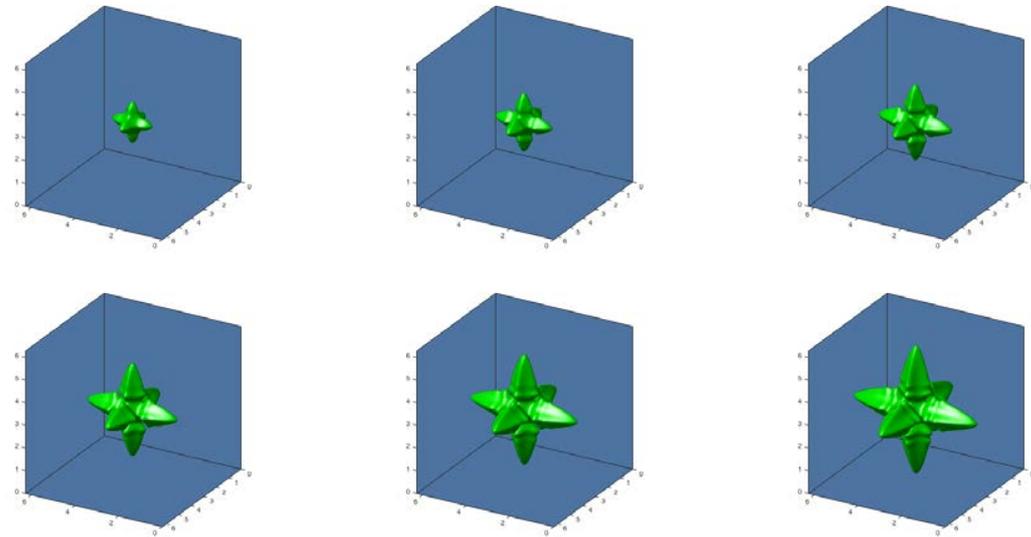

Figure 4.22. The dynamical evolution of dendritic crystal growth process in 3D with fourfold anisotropy and $K = 1.5$, computed by using the scheme SIEQ, the time step $\delta t = 1\mathrm{e}{-1}$, and parameters (4.7). Snapshots of the numerical approximation are taken at $T = 600$, 900, 1200, 1500, 1800 and 2100.

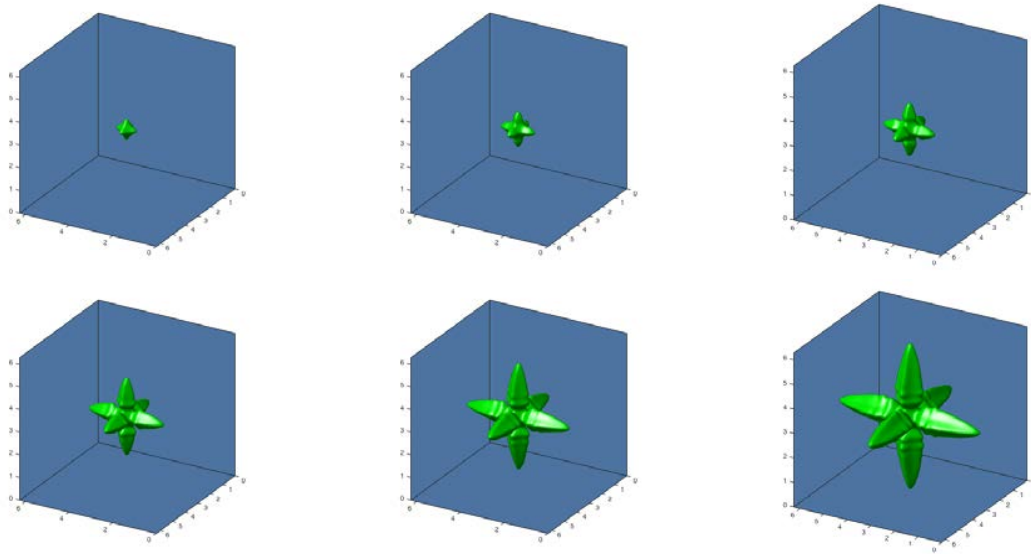

FIGURE 4.23. The dynamical evolution of dendritic crystal growth process in 3D with fourfold anisotropy and $K = 2$, computed by using the scheme SIEQ, the time step $\delta t = 1\text{e}-1$, and parameters (4.7). Snapshots of the numerical approximation are taken at $T = 300$, 600, 900, 1500, 2100, and 2700.

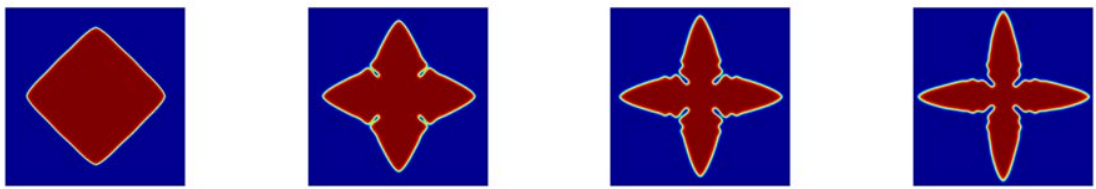

FIGURE 4.24. From left to right: 2D cross-section $\phi(\pi, \cdot, \cdot)$ for the computed solutions for $K = 0.5$, 1, 1.5, and 2. Snapshots are taken at $T = 900$, 1500, 2100, and 2700, respectively.

## 5. Concluding Remarks

In this paper, we have developed two efficient, semi-discrete in time, stabilized, second-order, linear schemes for solving the anisotropic dendritic Allen-Cahn phase field model. The first one is based on the linear stabilization approach where all nonlinear terms are treated explicitly and one only needs to solve two linear and decoupled second-order equations. The second one combines the recently developed Invariant Energy Quadratization approach with the linear stabilization approach. Two linear stabilization terms, which are shown to be crucial to remove the oscillations caused by the anisotropic coefficients numerically, are added as well. Compared to the existed schemes for the anisotropic model, the proposed SIEQ scheme that conquers the inconvenience from nonlinearities by linearizing the nonlinear terms in a new way, are provably unconditionally energy stable, and thus allow for large time steps in computations. We further numerically verify the accuracy in time and present various 2D and 3D numerical results for some benchmark numerical simulations.

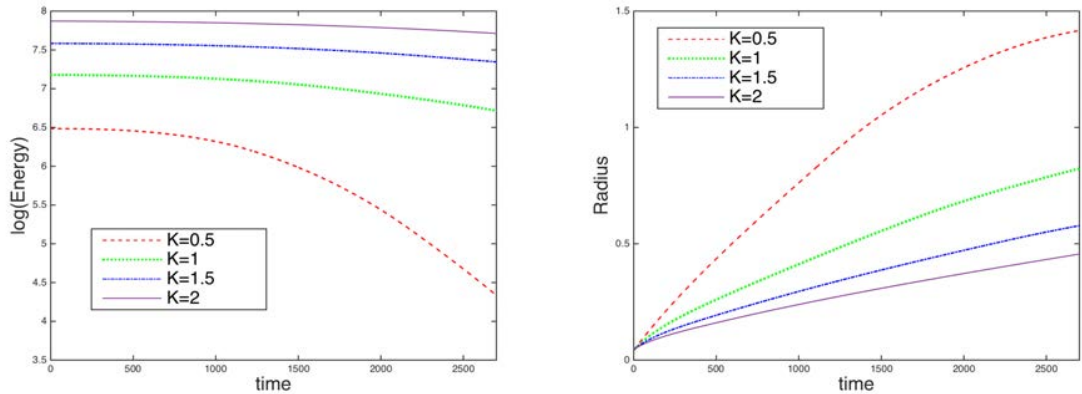

Figure 4.25. (a) Time evolutions of the logarithm of the free energy for the fourfold 3D example with $K = 0.5, 1, 1.5$, and $2$. (b) The radius of the crystals changing with time where the crystal radius is measured by an equivalent radius of a circle with the same area.

## References

[1] C. C. Chen, Y.L.Tsai, and C. W. Lan. Adaptive phase field simulation of dendritic crystal growth in a forced flow:2d vs. 3d morphologies. *Int. J. HeatMass Transfer*, 52:1158–1166, 2009.
[2] Q. Cheng, X. Yang, and J. Shen. Efficient and accurate numerical schemes for a hydro-dynamically coupled phase field diblock copolymer model. *J. Comp. Phys.*, 341:44–60, 2017.
[3] J. B. Collins and H. Levine. Diffuse interface model of diffusion-limited crystal growth. *Phys. Rev. B*, 31:6119, 1985.
[4] D. J. Eyre. Unconditionally gradient stable time marching the Cahn-Hilliard equation. In *Computational and mathematical models of microstructural evolution (San Francisco, CA, 1998)*, volume 529 of *Mater. Res. Soc. Sympos. Proc.*, pages 39–46. 1998.
[5] M. E. Glicksman, R. J. Schaefer, and J. D. Ayers. Dendritic growth-a test of theory. *Metallurgical Transactions A*, 7:1747–1759, 1976.
[6] B. I. Halperin, P. C. Hohenberg, and S.-K. Ma. Renormalization-group methods for critical dynamics: I. recursion relations and effects of energy conservation. *Phys. Rev. B*, 139:10, 1974.
[7] D. Han, A. Brylev, X. Yang, and Z. Tan. Numerical analysis of second order, fully discrete energy stable schemes for phase field models of two phase incompressible flows. *J. Sci. Comput.*, 70:965–989, 2017.
[8] Yasuji Sawada Haruo Honjo. Quantitative measurements on the morphology of a nh4br dendritic crystal growth in a capillary. *J. Crystal Growth*, 58:297–303, 1982.
[9] Y. He, Y. Liu, and T. Tang. On large time-stepping methods for the cahn-hilliard equation. *J. Appl. Num. Math.*, 57:616–628, 2007.
[10] S. C. Huang and M. E. Glicksman. Fundamentals of dendritic solidification ? i and ii. *Acta Metall*, 29:701?734, 1981.

[11] J.-H. Jeong, N. Goldenfeld, and J. A. Dantzig. Phase field model for three-dimensional dendritic growth with fluid flow. *Phys. Rev. E*, 55:041602, 2001.
[12] A. Karma and A. E. Lobkovsky. Unsteady crack motion and branching in a phase-field model of brittle fracture. *Phys. Rev. Lett.*, 92:245510, 2004.
[13] A. Karma and W. Rappel. Quantitative phase-field modeling of dendritic growth in two and three dimensions. *Phys. Rev. E*, 57:4323–4349, 1998.
[14] A. Karma and W. Rappel. Phase-field model of dendritic sidebranching with thermal noise. *Phys. Rev. E*, 60:3614–3625, 1999.
[15] R. Kobayashi. Modeling and numerical simulations of dendritic crystal growth. *Physica D*, 63:410, 1993.
[16] Y. Li and J. Kim. Phase-field simulations of crystal growth with adaptive mesh refinement. *Inter. J. Heat. Mass. Trans.*, 55:7926–7932, 2012.
[17] Y. Li, H. Lee, and J. Kim. A fast robust and accurate operator splitting method for phase field simulations of crystal growth. *Journal of Crystal Growth*, 321:176–182, 2011.
[18] E. Meca, V. Shenoy, and J. Lowengrub. Phase field modeling of two dimensional crystal growth with anisotropic diffusion. *Physical Review E*, 88:052409, 2013.
[19] B. Nestler, D. Danilov, and P. Galenko. Crystal growth of pure substances: phase field simulations in comparison with analytical and experimental results. *J. Comput. Phys.*, 207:221–239, 2005.
[20] M. Plapp and A. Karma. Multiscale finite-difference-diffusion-monte-carlo method for simulating dendritic solidification. *J. Comput. Phys*, 165:592–619, 2000.
[21] N. Provatas, N. Goldenfeld, and J. Dantzig. Efficient computation of dendritic microstructures using adaptive mesh refinement. *Phys. Rev. Lett*, 80:3308–3311, 1998.
[22] J. C. Ramirez, C. Beckermann, A. Karma, and H.-J. Diepers. Phase-field modeling of binary alloy solidication with coupled heat and solute diffusion. *Phys. Rev. E*, 69:051607, 2004.
[23] A. Shah, A. Haider, and S. Shah. Numerical simulation of two dimensional dendritic growth using phase field model. *World Journal of Mechanics*, 4:128–136, 2014.
[24] J. Shen, C. Wang, S. Wang, and X. Wang. Second-order convex splitting schemes for gradient flows with ehrlich-schwoebel type energy: application to thin film epitaxy. *SIAM J. Numer. Anal*, 50:105–125, 2012.
[25] J. Shen and J. Xu. Stabilized predictor-corrector schemes for gradient flows with strong anisotropic free energy. *Comm. Comput. Phys.*, to appear, 2018.
[26] J. Shen and X. Yang. Numerical Approximations of Allen-Cahn and Cahn-Hilliard Equations. *Disc. Conti. Dyn. Sys.-A*, 28:1669–1691, 2010.
[27] J. Shen and X. Yang. A phase-field model and its numerical approximation for two-phase incompressible flows with different densities and viscositites. *SIAM J. Sci. Comput.*, 32:1159–1179, 2010.
[28] S.-L. Wang, R. F. Sekerka, A. A. Wheeler, B. T. Murray, S.R. Coriell, R. J. Braun, and G.B. McFadden. Thermodynamically-consistent phase-field models for solidification. *Physica. D*, 69:189–200, 1993.
[29] J. A. Warren and W. J. Boettinger. Prediction of dentric growth and microsegregation patterns in a binary alloy using the phase field method. *Acta. metall. mater.*, 43:689–703, 1995.
[30] C. Xu and T. Tang. Stability analysis of large time-stepping methods for epitaxial growth models. *SIAM. J. Num. Anal.*, 44:1759–1779, 2006.
[31] X. Yang. Error analysis of stabilized semi-implicit method of allen-cahn equation. *Disc. Conti. Dyn. Sys.-B*, 11(4):1057–1070, 2009.
[32] X. Yang. Linear, first and second order and unconditionally energy stable numerical schemes for the phase field model of homopolymer blends. *J. Comput. Phys.*, 327:294–316, 2016.
[33] X. Yang. Numerical approximations for the cahn-hilliard phase field model of the binary fluid-surfactant system. *J. Sci. Comput*, 74:1533–1553, 2017.
[34] X. Yang and D. Han. Linearly first- and second-order, unconditionally energy stable schemes for the phase field crystal equation. *J. Comput. Phys.*, 330:1116–1134, 2017.

[35] X. Yang and L. Ju. Efficient linear schemes with unconditionally energy stability for the phase field elastic bending energy model. *Comput. Meth. Appl. Mech. Engrg.*, 315:691–712, 2017.

[36] X. Yang and L. Ju. Linear and unconditionally energy stable schemes for the binary fluid-surfactant phase field model. *Comput. Meth. Appl. Mech. Engrg.*, 318:1005–1029, 2017.

[37] X. Yang and G-D. Zhang. Numerical approximations of the Cahn-Hilliard and Allen-Cahn equations with general nonlinear potential using the invariant energy quadratization approach. *submitted*, 2018.

[38] X. Yang, J. Zhao, and Q. Wang. Numerical approximations for the molecular beam epitaxial growth model based on the invariant energy quadratization method. *J. Comput. Phys.*, 333:104–127, 2017.

[39] X. Yang, J. Zhao, Q. Wang, and J. Shen. Numerical approximations for a three components cahn-hilliard phase-field model based on the invariant energy quadratization method. *M3AS: Mathematical Models and Methods in Applied Sciences*, 27:1993–2030, 2017.

[40] J. Zhao, Q. Wang, and X. Yang. Numerical approximations for a phase field dendritic crystal growth model based on the invariant energy quadratization approach,. *Inter. J. Num. Meth. Engr.*, 110:279–300, 2017.